\pdfoutput=1  
\documentclass[11pt]{article}
\input{paper.sty}

\usepackage{amsmath,amssymb}
\usepackage{array}
\usepackage{bm}
\usepackage{booktabs}
\usepackage[section]{placeins}
\usepackage[title]{appendix}
\usepackage{graphicx}
\usepackage{tikz}
\usetikzlibrary{calc,arrows.meta,positioning,decorations.pathreplacing,shapes.geometric,backgrounds,patterns,3d,fit}
\usepackage[hidelinks]{hyperref}
\usepackage[capitalize]{cleveref}
\usepackage{microtype}

\newcommand{\kron}{\otimes}
\newcommand{\Dx}{D_{x}}
\newcommand{\Dy}{D_{y}}
\newcommand{\Dz}{D_{z}}
\newcommand{\Dxx}{D_{xx}}
\newcommand{\Dyy}{D_{yy}}
\newcommand{\Dzz}{D_{zz}}
\newcommand{\Ix}{I_{x}}
\newcommand{\Iy}{I_{y}}
\newcommand{\Iz}{I_{z}}
\newcommand{\Ax}{A_{x}}
\newcommand{\Ay}{A_{y}}
\newcommand{\Az}{A_{z}}
\newcommand{\Axx}{A_{xx}}
\newcommand{\Rx}{R_{x}}
\newcommand{\Ry}{R_{y}}
\newcommand{\Rz}{R_{z}}
\newcommand{\Rxx}{R_{xx}}
\newcommand{\R}{\mathbb{R}}
\newcommand{\vect}{\operatorname{vec}}
\newcommand{\diag}{\operatorname{diag}}
\newcommand{\dt}{\Delta t}
\newcommand{\bigO}{\mathcal{O}}
\newcommand{\Lap}{\Delta_h}
\newcommand{\Nall}{N}

\title{No 3D Matrices:\\
A Unified Tensor-Product View of Matrix-Free Cartesian PDE Solvers}

\author{%
  Yong Yi Bay\thanks{Equal contribution. Emails: \texttt{\{yongyibay,kallie.a.yearick\}@gmail.com}.}
  \hspace{2.5em}
  Kathleen A. Yearick\footnotemark[1] \\
  {\normalfont\normalsize PhD, University of Illinois at Urbana-Champaign}
}

\begin{document}
\raggedbottom
\maketitle

\begin{abstract}
Every Cartesian three-dimensional PDE solver hides a structural
secret that production CFD codes have used for half a century and
that graduate-level textbooks rarely state plainly.  The derivative
matrices, the compact Pad\'{e} line solves, the Galerkin mass
inversions, the alternating-direction-implicit substeps, and even the
fast Poisson and Helmholtz diagonalization transforms all factor
along the coordinate axes and collapse into repeated one-dimensional
banded kernels executed along the grid lines.  The three-dimensional
operator exists only on paper; it is never assembled, factored, or
stored.  This paper is the manual for that collapse.  We derive the
Kronecker-product algebra that makes it exact, carry it cleanly
through central differences, compact schemes, tensor-product
Galerkin, B-spline and isogeometric methods, collocation, ADI time
stepping, and direct Poisson and Helmholtz solves, and bring into the
open the three production tricks that turn the reduction into
hardware-conscious floating-point throughput on real machines: the
multi-right-hand-side reshape that exposes a sweep as one batched
line kernel (a dense BLAS-3 GEMM when the line factor is dense or
element-local, a banded or stencil kernel when it is not), the sum
factorization that rescues high-order Galerkin from the
$\bigO(p^{2d})$ quadrature trap, and the pencil decomposition that
keeps every direction contiguous across an MPI cluster.  For fixed
stencil width or fixed polynomial degree, the compute cost stays
$\bigO(\Nall)$ in the total number of unknowns
$\Nall = N_x N_y N_z$; the operator storage drops to
$\bigO(N_x + N_y + N_z)$ up to bandwidth constants; direct separable
Poisson and Helmholtz solvers add the expected transform cost; the
line kernels are embarrassingly parallel.  These facts are familiar
to practitioners but rarely assembled in one place; this paper
collects them and shows how to use them.
\end{abstract}


\section{Introduction}
\label{sec:intro}

Imagine a student sitting down to solve the heat equation on a box
with $200$ points per side.  The Laplacian inside an implicit time
step would be an $8{,}000{,}000 \times 8{,}000{,}000$ matrix if it
were written monolithically.  Dense storage would require roughly
$5.1\times10^{14}$ bytes, about half a petabyte in double precision.
Even the seven-point sparse matrix has about $5.6\times10^7$
nonzeros, already close to a gigabyte before any solver fill-in or
preconditioner storage, and a direct factorization of that object
would still occupy a serious workstation for the better part of a
night.  This sounds like a hard problem.  It is not.

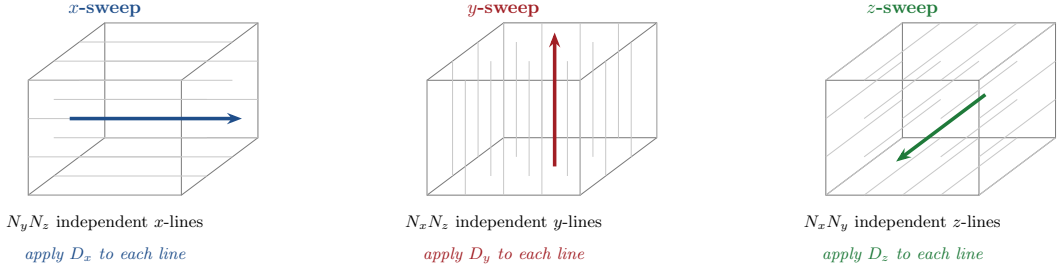
\begin{figure}[t]
  \centering
  \resizebox{\linewidth}{!}{%
  \begin{tikzpicture}[
    x={(0.72cm,0cm)}, y={(0cm,0.72cm)}, z={(-0.48cm,-0.36cm)},
    every node/.style={font=\small},
  ]
    \pgfmathsetmacro{\Nx}{4}
    \pgfmathsetmacro{\Ny}{3}
    \pgfmathsetmacro{\Nz}{3}
    \definecolor{sweepblue}{HTML}{1f4e8a}
    \definecolor{sweepred}{HTML}{a11e24}
    \definecolor{sweepgreen}{HTML}{1f7a3a}
    \definecolor{boxedge}{HTML}{7a7a7a}
    \definecolor{faintline}{HTML}{c8c8c8}

    \newcommand{\wireframe}{%
      \draw[boxedge, line width=0.5pt] (0,0,0) -- (\Nx,0,0) -- (\Nx,\Ny,0) -- (0,\Ny,0) -- cycle;
      \draw[boxedge, line width=0.5pt] (0,0,\Nz) -- (\Nx,0,\Nz) -- (\Nx,\Ny,\Nz) -- (0,\Ny,\Nz) -- cycle;
      \draw[boxedge, line width=0.5pt] (0,0,0) -- (0,0,\Nz);
      \draw[boxedge, line width=0.5pt] (\Nx,0,0) -- (\Nx,0,\Nz);
      \draw[boxedge, line width=0.5pt] (\Nx,\Ny,0) -- (\Nx,\Ny,\Nz);
      \draw[boxedge, line width=0.5pt] (0,\Ny,0) -- (0,\Ny,\Nz);
    }

    \begin{scope}[xshift=0cm]
      \wireframe
      \foreach \k in {0,...,\Nz}{
        \foreach \j in {0,...,\Ny}{
          \draw[faintline, line width=0.45pt] (0,\j,\k) -- (\Nx,\j,\k);
        }
      }
      \draw[sweepblue, line width=1.9pt, -{Stealth[length=2.8mm]}]
        (-0.25,1,1) -- (\Nx+0.25,1,1);
      \node[anchor=south, sweepblue, font=\bfseries]
        at (\Nx/2,\Ny+1.5,\Nz) {$x$-sweep};
      \node[anchor=north] at (\Nx/2,-0.35,\Nz) {$N_y N_z$ independent $x$-lines};
      \node[anchor=north, sweepblue, font=\small\itshape]
        at (\Nx/2,-1.2,\Nz) {apply $D_x$ to each line};
    \end{scope}

    \begin{scope}[xshift=7.5cm]
      \wireframe
      \foreach \k in {0,...,\Nz}{
        \foreach \i in {0,...,\Nx}{
          \draw[faintline, line width=0.45pt] (\i,0,\k) -- (\i,\Ny,\k);
        }
      }
      \draw[sweepred, line width=1.9pt, -{Stealth[length=2.8mm]}]
        (2,-0.25,1) -- (2,\Ny+0.25,1);
      \node[anchor=south, sweepred, font=\bfseries]
        at (\Nx/2,\Ny+1.5,\Nz) {$y$-sweep};
      \node[anchor=north] at (\Nx/2,-0.35,\Nz) {$N_x N_z$ independent $y$-lines};
      \node[anchor=north, sweepred, font=\small\itshape]
        at (\Nx/2,-1.2,\Nz) {apply $D_y$ to each line};
    \end{scope}

    \begin{scope}[xshift=15.0cm]
      \wireframe
      \foreach \j in {0,...,\Ny}{
        \foreach \i in {0,...,\Nx}{
          \draw[faintline, line width=0.45pt] (\i,\j,0) -- (\i,\j,\Nz);
        }
      }
      \draw[sweepgreen, line width=1.9pt, -{Stealth[length=2.8mm]}]
        (2,1,-0.25) -- (2,1,\Nz+0.25);
      \node[anchor=south, sweepgreen, font=\bfseries]
        at (\Nx/2,\Ny+1.5,\Nz) {$z$-sweep};
      \node[anchor=north] at (\Nx/2,-0.35,\Nz) {$N_x N_y$ independent $z$-lines};
      \node[anchor=north, sweepgreen, font=\small\itshape]
        at (\Nx/2,-1.2,\Nz) {apply $D_z$ to each line};
    \end{scope}
  \end{tikzpicture}}
  \caption{The whole paper in one picture.  A Cartesian three-dimensional
  operator is executed as a collection of independent one-dimensional
  line operations.  For an $x$-derivative, the same small banded matrix
  $D_x$ is applied to each $x$-line of the grid, giving $N_y N_z$
  independent line kernels (left).  The $y$- and $z$-derivatives reuse
  the same rule with the direction relabeled (center, right).  No
  three-dimensional matrix is ever assembled: the kernel is
  one-dimensional, and the dimension only counts how many lines that
  kernel must visit.}
  \label{fig:hero-sweep}
\end{figure}

Production Cartesian solvers do not build that object.  They exploit
one identity, which a few pages of Kronecker algebra will make exact.
With $D_{xx}$, $D_{yy}$, and $D_{zz}$ the three one-dimensional
second-derivative matrices (each of size $200$, tridiagonal, a few
kilobytes total),
\begin{equation*}
  \Lap
  \;=\;
  I_z \kron I_y \kron D_{xx}
  \;+\;
  I_z \kron D_{yy} \kron I_x
  \;+\;
  D_{zz} \kron I_y \kron I_x.
\end{equation*}
The $8{,}000{,}000\times 8{,}000{,}000$ matrix on the left is the
sum of three Kronecker products of the three small matrices on the
right.  We never build the object on the left.  Applying $\Lap$ to a
field is three sweeps along the three families of grid lines, each
sweep a batch of $200$-entry tridiagonal matvecs.  A factored
implicit method, such as ADI, replaces the coupled solve by the
corresponding sequence of one-dimensional banded solves; a separable
Poisson or Helmholtz solve uses the same directional structure through
fast diagonalization.  The algorithm lives inside the
one-dimensional pieces.  This is the picture behind every serious
Cartesian solver, and the one \cref{fig:hero-sweep}
draws for the eye.

The idea is old in the best possible way.  It has lived in
direct-numerical-simulation
solvers~\citep{bewley2001dns,kim1985application}, in FFT
Poisson packages~\citep{costa2018fft}, and in spectral-element
infrastructure for incompressible flow and
acoustics~\citep{deville2002high} for decades.  It is also the design
logic behind pencil-decomposition libraries such as
2DECOMP\&FFT~\citep{li2010decomp,rolfo2023decomp}.  What has been
missing is a single place where the argument is written down cleanly,
once, for every discretization the reader is likely to care about,
where the scope of the reduction is stated honestly, and where the
production tricks that turn the argument into running floating-point
code are ushered into the open.

The bookkeeping that makes the picture exact is the Kronecker
product, which Van Loan memorably called
``ubiquitous''~\citep{loan2000ubiquitous} and which Strang uses as
the running thread of computational linear
algebra~\citep{strang2007computational}.  Two identities do most of
the work: the mixed-product rule $(A \kron B)(C \kron D) = AC \kron BD$,
and the inverse rule $(A \kron B)^{-1} = A^{-1} \kron B^{-1}$.  The
first composes directional operators without ever forming their
product.  The second explains why a tensor-product mass matrix or
compact line factor can be inverted by independent one-dimensional
solves.  Every subsequent section of this paper is a variation on
those two identities plus the column-major vectorization rule.

The counting that makes the picture \emph{fast} is even simpler.  A
one-dimensional banded matrix-vector product with bandwidth $w_x$
costs $\bigO(w_xN_x)$.  A three-dimensional $x$-sweep is $N_yN_z$
such products, hence $\bigO(w_xN_xN_yN_z)=\bigO(w_x\Nall)$.  For
fixed stencil width or fixed polynomial degree this is simply
$\bigO(\Nall)$.  Banded line solves have the same linear count, up to
bandwidth constants, and fast transforms add their usual logarithmic
factor.  The full operator is large only if it is assembled; the line
kernels are never large.

\begin{table}[t]
\centering
\footnotesize
\renewcommand{\arraystretch}{1.2}
\begin{tabular}{@{}p{0.25\textwidth}p{0.34\textwidth}p{0.28\textwidth}@{}}
\toprule
\textbf{Task} & \textbf{Tensor-product form} & \textbf{How it runs} \\
\midrule
Derivative sweep & $I_z\kron I_y\kron D_x$ and relabelings & independent banded line applies, $\bigO(\Nall)$ for fixed bandwidth \\
\addlinespace[3pt]
Compact derivative & $(I_z\kron I_y\kron A_x)v=(I_z\kron I_y\kron R_x)u$ & banded right-hand side plus banded line solve \\
\addlinespace[3pt]
Tensor-product mass inverse & $(M_z\kron M_y\kron M_x)^{-1}$ & three families of 1D mass solves \\
\addlinespace[3pt]
ADI substep & $I-r(I_z\kron I_y\kron D_{xx})$ & one directional batch of line solves \\
\addlinespace[3pt]
Fast Poisson solve & diagonalized Kronecker sum & transforms by direction, scalar division, inverse transforms \\
\bottomrule
\end{tabular}
\caption{The recurring pattern.  The algebra names the operator; the
implementation visits one-dimensional lines.  The field storage is
still $\bigO(\Nall)$, but the \emph{operator} storage is only the set
of one-dimensional factors, typically $\bigO(N_x+N_y+N_z)$ up to
bandwidth constants.}
\label{tab:pattern}
\end{table}

That combination, a Kronecker factorization for the algebra and a
linewise count for the arithmetic, runs the whole paper.
\Cref{sec:example} makes both visible on a $4 \times 3$ grid small
enough to hold on one page, and recasts the 2D Poisson equation as a
Sylvester matrix equation that makes the line action plain.
\Cref{sec:kron} collects the two identities and fixes notation.
\Cref{sec:1d-ops,sec:3d-ops} build the operator family from the one
dimensional pieces: uniform and stretched grids, derivatives of all
orders, mixed derivatives, the Laplacian.  \Cref{sec:compact}
handles compact Pad\'{e} schemes and their implicit line
structure~\citep{lele1992compact,beam1976implicit}.
\Cref{sec:galerkin} extends the algebra to tensor-product Galerkin,
B-spline, isogeometric, and collocation
methods~\citep{deboor1978practical,hughes2005isogeometric,cottrell2009isogeometric,johnson2005higher,bay2019dissertation}.
\Cref{sec:implicit} folds in implicit time stepping through the
classical Peaceman--Rachford and Douglas--Gunn ADI
splittings~\citep{peaceman1955numerical,douglas1962alternating,douglas1964numerical}.
\Cref{sec:fast-diag} gives the direct-solver route by the
Lynch--Rice--Thomas fast-diagonalization
construction~\citep{lynch1964direct,haidvogel1979accurate}.
\Cref{sec:applicability} marks the boundary of what the framework
can reach.  \Cref{sec:practice} sets out the three
implementation tricks that turn the algebra into efficient,
hardware-conscious kernels on real machines.  \Cref{sec:experiment}
measures the gap between assembling the three-dimensional matrix and
never forming it on a model 3D Poisson solve.  \Cref{sec:summary} is the
implementation recipe in one page.  The appendix collects the short
derivations and the $\vect{}$ identity that translate Kronecker
formulas into explicit loops.

\section{A two-dimensional motivating example}
\label{sec:example}

Small cases reveal the pattern.  Let us take a
2D grid with $N_x = 4$ points in $x$ and $N_y = 3$ points in $y$
(boundary details ignored for clarity).  The grid has $4 \times 3 = 12$
points total, which is small enough to write every matrix on one page.

In 1D, the second-order central-difference first derivative on 4 points
with spacing $h$ is the tridiagonal matrix
\begin{equation}
\label{eq:Dx-small}
  \Dx = \frac{1}{2h}
  \begin{bmatrix}
     0 &  1 &  0 &  0 \\
    -1 &  0 &  1 &  0 \\
     0 & -1 &  0 &  1 \\
     0 &  0 & -1 &  0
  \end{bmatrix}
  \in \R^{4 \times 4}.
\end{equation}
(One-sided stencils at the first and last rows would close the
boundaries in a real code.  The specific closure does not affect
anything that follows.)

On the 2D grid, the objective is to compute $\partial u / \partial x$
at all 12 points.  Because the grid is Cartesian, the $x$-stencil at
point $(i,j)$ depends only on the $x$-neighbors $(i-1,j)$ and
$(i+1,j)$, at the same $y$-index $j$.  For each $j$, $\Dx$ is therefore
applied to the four values along the corresponding $x$-line:

\medskip
\begin{tabular}{ll}
  Line $j=1$: & $\Dx\begin{bmatrix} u_{1,1} \\ u_{2,1} \\ u_{3,1} \\ u_{4,1} \end{bmatrix}$ \\[12pt]
  Line $j=2$: & $\Dx\begin{bmatrix} u_{1,2} \\ u_{2,2} \\ u_{3,2} \\ u_{4,2} \end{bmatrix}$ \\[12pt]
  Line $j=3$: & $\Dx\begin{bmatrix} u_{1,3} \\ u_{2,3} \\ u_{3,3} \\ u_{4,3} \end{bmatrix}$
\end{tabular}
\medskip

The operation therefore consists of three independent applications of
the same $4 \times 4$ matrix.  No two-dimensional matrix is formed.

\subsection{Three dimensions is only a relabeling}
\label{sec:pseudocode}

Moving the same argument to 3D changes nothing conceptual.  The
$x$-line, the $y$-line, and the $z$-line each play the role that the
single $x$-line played above; each of them is swept by its own
tridiagonal; each sweep is then stitched into the full field.
\Cref{fig:pseudocode} is a one-page template for the whole family.
Every section that follows instantiates these two building blocks,
\texttt{sweep} and \texttt{solve}, with a different discretization.

\begin{figure}[t]
\hrule\medskip
\small
\begin{verbatim}
  # --- Building block: apply 1D operator along one direction ---

  function sweep_x(D_x, u):            # explicit: multiply
      for k = 1, ..., N_z:
          for j = 1, ..., N_y:
              v(:,j,k) = D_x @ u(:,j,k)      # 1D matvec, size N_x
      return v

  function solve_x(A_x, R_x, u):       # implicit: banded line solve
      for k = 1, ..., N_z:
          for j = 1, ..., N_y:
              rhs       = R_x @ u(:,j,k)      # 1D matvec
              v(:,j,k)  = banded_solve(A_x, rhs)
      return v

  # sweep_y and sweep_z are identical, looping over the other two indices.

  # --- Derivatives ---

  du_dx  = sweep_x(D_x, u)             # first derivative in x
  du_dy  = sweep_y(D_y, u)             # first derivative in y
  d2u_dx2 = sweep_x(D_xx, u)           # second derivative in x

  # --- Laplacian: three sweeps, then add ---

  lap_u = sweep_x(D_xx, u) + sweep_y(D_yy, u) + sweep_z(D_zz, u)

  # --- Compact (Pade) derivative: same loop, implicit version ---

  du_dx = solve_x(A_x, R_x, u)         # A_x * du_dx = R_x * u, per line

  # --- Mixed derivative d^2u/dxdy: two sweeps in sequence ---

  tmp    = sweep_x(D_x, u)             # differentiate in x
  d2u_dxdy = sweep_y(D_y, tmp)         # then differentiate in y
\end{verbatim}
\medskip\hrule
\caption{Tensor-product line-sweep template.  Every 3D operation is
  expressed as a loop of 1D operations along grid lines.  The explicit
  variant
  (\texttt{sweep}) is a matrix-vector product; the implicit variant
  (\texttt{solve}) is a banded line solve.  The Laplacian is three
  sweeps summed.  ADI time stepping (\cref{sec:implicit}) and fast
  diagonalization (\cref{sec:fast-diag}) compose these same building
  blocks.}
\label{fig:pseudocode}
\end{figure}

If this action is written as a single matrix acting on all 12 unknowns
at once (stacking the lines so that $x$ varies fastest), the matrix
would be the $12 \times 12$ block-diagonal matrix
\begin{equation}
\label{eq:block-diag}
  \begin{bmatrix}
    \Dx &     &     \\
        & \Dx &     \\
        &     & \Dx
  \end{bmatrix}
  = I_3 \kron \Dx,
\end{equation}
where $I_3$ is the $3 \times 3$ identity (one row per $y$-line) and
$\kron$ denotes the Kronecker product.  The notation
``$I_3 \kron \Dx$'' is a compact representation of ``apply $\Dx$ to
each line independently.''  This linewise action is the core of the
framework.

\paragraph{Poisson is a Sylvester equation in disguise.}
There is a second picture of the same algebra, and it is the one
mathematicians will recognize first.  Arrange the unknowns into a
rectangular array $U \in \R^{N_x \times N_y}$, one column per
$y$-line, and do not vectorize at all.  The Laplacian acts by
columns and by rows,
\begin{equation}
\label{eq:sylvester-view}
  \Lap\, U \;=\; \Dxx\, U \;+\; U\, \Dyy^T,
\end{equation}
because $\Dxx$ differentiates down the columns and $\Dyy^T$
differentiates across the rows.  The 2D Poisson problem
$-\Lap u = f$ is therefore the classical \emph{Sylvester matrix
equation}
\begin{equation}
\label{eq:sylvester}
  -\Dxx\, U - U\, \Dyy^T \;=\; F,
\end{equation}
and the storage collapses to two small matrices and a rectangular
grid of data: nothing of size $\Nall \times \Nall$ appears anywhere.
The equivalent long-vector form $-(\Iy \kron \Dxx + \Dyy \kron
\Ix)\,\vect(U) = \vect(F)$ says the same thing
(\cref{app:vec}); the rectangular picture simply makes the directional
action visible.  Diagonalizing $\Dxx$ and $\Dyy$ in their own
eigenbases uncouples the whole equation pointwise; that is the
fast-diagonalization route of
\cref{sec:fast-diag}~\citep{lynch1964direct}.  The picture
generalizes to three dimensions without change: the 3D Poisson
problem is a three-term matrix equation on a tensor, with $\Dxx$,
$\Dyy$, and $\Dzz$ acting each along its own coordinate axis.

\section{Kronecker-product identities used throughout}
\label{sec:kron}

The Kronecker product $A \kron B$ of an $m \times n$ matrix $A$ and a
$p \times q$ matrix $B$ is the $mp \times nq$ block matrix formed by
replacing each entry $a_{ij}$ of $A$ with the block $a_{ij}\,B$:
\begin{equation}
\label{eq:kron-def}
  A \kron B
  =
  \begin{bmatrix}
    a_{11} B & a_{12} B & \cdots & a_{1n} B \\
    a_{21} B & a_{22} B & \cdots & a_{2n} B \\
    \vdots   & \vdots   & \ddots & \vdots   \\
    a_{m1} B & a_{m2} B & \cdots & a_{mn} B
  \end{bmatrix}.
\end{equation}
Equation~\eqref{eq:block-diag} shows that $I_3 \kron \Dx$ produces the
same block-diagonal matrix: the identity has ones on its diagonal, so
each diagonal block is $1 \cdot \Dx = \Dx$, and the off-diagonal blocks
are $0 \cdot \Dx = 0$.

Only two standard identities are required for the developments that
follow.

\paragraph{Property 1: the mixed-product rule.}
If the matrix products $AC$ and $BD$ are defined, then
\begin{equation}
\label{eq:mixed-product}
  (A \kron B)(C \kron D) = (AC) \kron (BD).
\end{equation}
That is, the corresponding factors multiply independently.  This
identity underlies the composition of directional operators.

\paragraph{Property 2: the inverse rule.}
If $A$ and $B$ are both invertible, then
\begin{equation}
\label{eq:kron-inv}
  (A \kron B)^{-1} = A^{-1} \kron B^{-1}.
\end{equation}
This follows from the mixed-product rule:
$(A \kron B)(A^{-1} \kron B^{-1}) = (AA^{-1}) \kron (BB^{-1}) = I \kron I = I$.

These identities, together with the vectorization convention below,
are enough for every construction in the paper.  They belong to a
much larger body of Kronecker-product algebra treated systematically
by \citet{golub2013matrix,loan2000ubiquitous} and, in the
higher-order-tensor setting, by \citet{kolda2009tensor}.
Appendix~\ref{app:linalg} records the short derivations and the
$\vect{}$ identity that converts the abstract Kronecker formulas into
explicit line sweeps.

\paragraph{Vectorization convention.}
One bookkeeping detail concerns how a 2D array $U(i,j)$ is flattened
into a vector $\bm{u}$: columns are stacked so that the $x$-index
varies fastest.
For a $4 \times 3$ grid:
\begin{equation}
  \bm{u} = \bigl[
    u_{1,1},\; u_{2,1},\; u_{3,1},\; u_{4,1},\;
    u_{1,2},\; u_{2,2},\; \ldots,\;
    u_{4,3}
  \bigr]^T.
\end{equation}
In 3D, the ordering is $x$ fastest, then $y$, then $z$.  This
column-major ordering is what makes the Kronecker product formulas come
out with the rightmost factor acting on $x$, the middle factor on $y$,
and the leftmost on $z$.

\section{One-dimensional finite-difference operators}
\label{sec:1d-ops}

\begin{figure}[t]
  \centering
  \includegraphics[width=\linewidth]{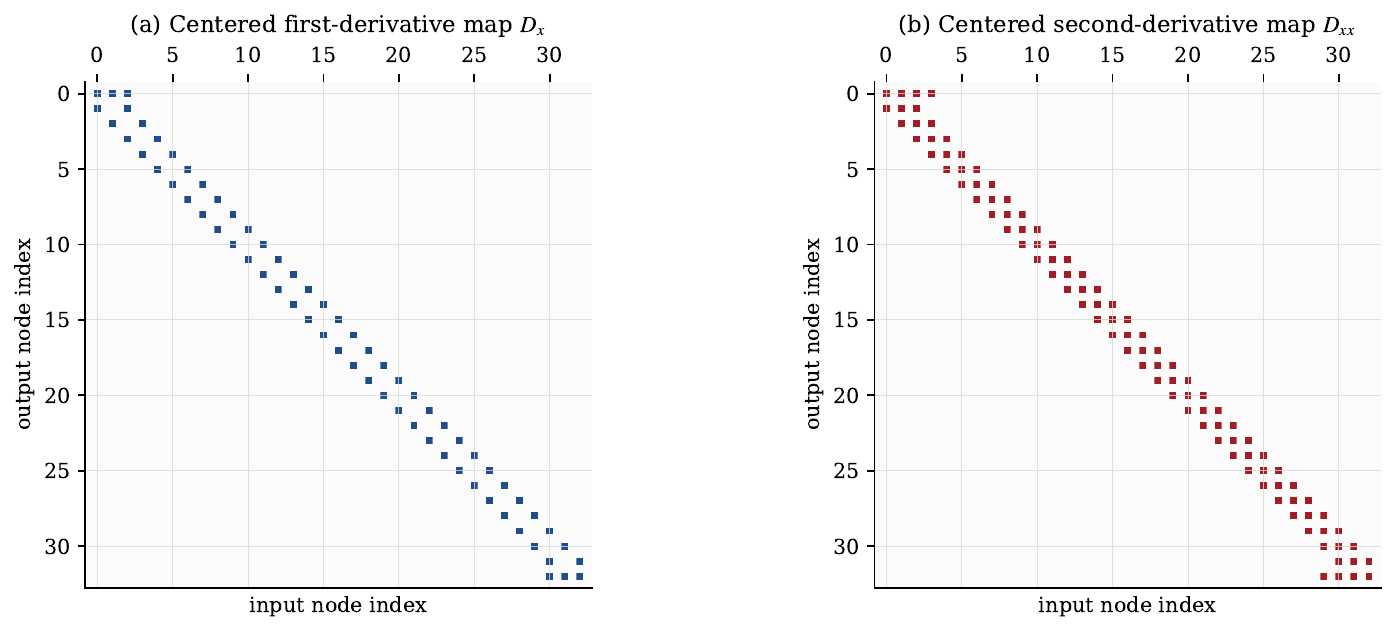}
  \caption{Centered finite-difference maps in 1D.  The first- and
  second-derivative operators remain banded, so a line apply touches
  each unknown a constant number of times.  This is the base case for
  the tensor-product reduction: local operators remain local in each
  coordinate direction.}
  \label{fig:fd-sparsity}
\end{figure}

This section collects the 1D building blocks.  On a non-uniform grid,
the operator remains banded and only the entries change.
\Cref{fig:fd-sparsity} emphasizes this structure before the formulas
are introduced.  The first- and second-derivative maps occupy thin
bands around the diagonal.  Consequently, a full 3D sweep remains
linear in $\Nall$: every output value depends on only a fixed number of
nearby inputs, with no long-range fill-in and no dense intermediate
state.

\subsection{Uniform grids}
\label{sec:1d-uniform}

On a uniform grid $x_i = x_0 + i\,h$, the standard central-difference
approximations are:
\begin{align}
  \text{First derivative:} \quad
  u'(x_i) &\approx \frac{u_{i+1} - u_{i-1}}{2h},
  \label{eq:fd1} \\[4pt]
  \text{Second derivative:} \quad
  u''(x_i) &\approx \frac{u_{i-1} - 2u_i + u_{i+1}}{h^2}.
  \label{eq:fd2}
\end{align}
Collecting these for all grid points gives the tridiagonal matrices
$\Dx$ and $\Dxx$.  Applying either to a vector of length $N_x$ costs
$\bigO(N_x)$ operations.  Higher-order stencils (fourth-order,
sixth-order) widen the bandwidth but do not change the
structure~\citep{fornberg1988generation}.  In three dimensions, the
operator is obtained by repeating the same one-dimensional banded apply
on many parallel lines.

\paragraph{Boundary closures.}
Equations~\eqref{eq:fd1}--\eqref{eq:fd2} apply at interior points.
The endpoint rows of $\Dx$ and $\Dxx$ require boundary closures.  A
standard second-order first-derivative closure uses
$u'(x_1) \approx (-3 u_1 + 4 u_2 - u_3)/(2h)$ and the reflected
formula at $x_{N_x}$; a second-order second-derivative closure may use
$(2u_1-5u_2+4u_3-u_4)/h^2$ at the left endpoint and the reflected
row at the right endpoint.  Compact, Galerkin, and B-spline
discretizations carry their own closure families
(\citealp{lele1992compact}; \citealp{deboor1978practical}).  The
choice of closure modifies only a small number of boundary rows in
the 1D operator, leaves the interior banded structure intact, and
inherits into the 3D apply with no change to the Kronecker
factorization.

\paragraph{The line solver: Thomas algorithm.}
Every implicit line operation in this paper invokes a banded triangular
solve, and for tridiagonal systems the workhorse is the Thomas
algorithm~\citep{thomas1949elliptic}.  Given $\Ax\bm{x} = \bm{d}$
with $\Ax$ tridiagonal of subdiagonal $a_i$, diagonal $b_i$, and
superdiagonal $c_i$, the forward sweep eliminates the subdiagonal.
With $\tilde c_1=c_1/b_1$ and $\tilde d_1=d_1/b_1$, the interior rows
$i=2,\ldots,N_\xi-1$ use
\begin{equation}
\label{eq:thomas-fwd}
  m_i = b_i-a_i\tilde c_{i-1},\quad
  \tilde c_i = \frac{c_i}{m_i},\quad
  \tilde d_i = \frac{d_i-a_i\tilde d_{i-1}}{m_i},
\end{equation}
and the final row closes the sweep with
$m_{N_\xi} = b_{N_\xi}-a_{N_\xi}\tilde c_{N_\xi-1}$ and
$\tilde d_{N_\xi} = (d_{N_\xi}-a_{N_\xi}\tilde d_{N_\xi-1})/m_{N_\xi}$.
These formulas assume the pivots $m_i$ are nonzero.  That is the
usual case for the diagonally dominant elliptic line factors that
motivate the discussion; otherwise one must use pivoting or a more
general banded solver.
The backward sweep recovers the solution:
\begin{equation}
\label{eq:thomas-bwd}
  x_{N_\xi} = \tilde d_{N_\xi},
  \qquad
  x_i = \tilde d_i - \tilde c_i\, x_{i+1},
  \qquad
  i = N_\xi-1, \ldots, 1.
\end{equation}
The total cost is $\bigO(N_\xi)$ operations (a small constant number
of flops per row) and only $\bigO(N_\xi)$ memory traffic per line.  A
banded solve with semibandwidth $w$ replaces the two scalar updates
above with $\bigO(w^2)$ local elimination work per row for a general
banded factorization, or $\bigO(w)$ work when a fixed-band LU factor
has already been computed.  In either case, for fixed $w$ the line
solve is linear in $N_\xi$.  A 3D implicit sweep is
$\Nall/N_\xi$ independent invocations of this two-pass kernel.

\subsection{Non-uniform grids}
\label{sec:1d-nonuniform}

Now let the spacing vary: $h_i^{+} = x_{i+1} - x_i$ and
$h_i^{-} = x_i - x_{i-1}$.  The three-point central formulas extend
to
\begin{align}
  u'(x_i) &\approx
  \frac{
    -(h_i^{+})^2 \, u_{i-1}
    + \bigl[(h_i^{+})^2 - (h_i^{-})^2\bigr] u_i
    + (h_i^{-})^2 \, u_{i+1}
  }{
    h_i^{+}\, h_i^{-}\, (h_i^{+} + h_i^{-})
  },
  \label{eq:fd1-nonunif} \\[4pt]
  u''(x_i) &\approx
  \frac{2}{h_i^{+} + h_i^{-}}
  \left(
    \frac{u_{i+1} - u_i}{h_i^{+}} - \frac{u_i - u_{i-1}}{h_i^{-}}
  \right).
  \label{eq:fd2-nonunif}
\end{align}
On a uniform grid, \eqref{eq:fd1-nonunif} reduces to
$(u_{i+1} - u_{i-1})/(2h)$ and \eqref{eq:fd2-nonunif} reduces to
$(u_{i-1} - 2 u_i + u_{i+1})/h^2$.  A Taylor expansion gives two
useful cautions.  The first-derivative formula
\eqref{eq:fd1-nonunif} is exact for quadratics and has a second-order
local truncation error when the adjacent spacings are comparable.  The
three-point second-derivative formula \eqref{eq:fd2-nonunif} is also
exact for quadratics, but its leading cubic error is proportional to
$h_i^{+}-h_i^{-}$.  Thus it is second order on smoothly mapped grids
where $h_i^{+}-h_i^{-}=\bigO(h^2)$, but only first order on a
strongly irregular grid where that difference is $\bigO(h)$.  If one
needs second-order accuracy on arbitrary non-uniform point sets, a
wider stencil generated by Fornberg's algorithm~\citep{fornberg1988generation}
or a smooth coordinate transform is the safer choice.  In every case
the stencil weights vary from point to point, but the operator remains
banded.

For the purposes of this paper, the essential point is that a
non-uniform Cartesian grid is still a tensor product of 1D grids, so
the Kronecker product structure is unchanged.  The one-dimensional
matrices acquire non-constant diagonals, but the geometry may stretch,
cluster, or bias the points without changing the cost model: banded 1D
operators still give linear-time sweeps in the total number of grid
values.

\section{From 1D to 3D: every operator, one identity at a time}
\label{sec:3d-ops}

The same Kronecker machinery produces the whole operator family.
Throughout this section, $\Ix$, $\Iy$, $\Iz$ are identity matrices
of sizes $N_x$, $N_y$, $N_z$, and $\bm{u} \in \R^{\Nall}$ is the
flattened field with $x$ varying fastest.  The rule for reading any
triple Kronecker product is the same in every case: the rightmost
factor acts on the fastest ($x$) index, the middle factor on $y$,
the leftmost on $z$.

\subsection{First derivatives}
\label{sec:3d-first}

The central-difference approximation to $\partial u/\partial x$ at
the interior point $(i,j,k)$,
\begin{equation}
\label{eq:3d-dx-pointwise}
  \left(\frac{\partial u}{\partial x}\right)_{ijk}
  \approx
  \frac{u_{i+1,j,k} - u_{i-1,j,k}}{2h_x},
\end{equation}
touches only the two $x$-neighbors at the same $(j,k)$.  The
directional Kronecker placement is therefore dictated by
\cref{eq:3d-dx-pointwise} itself:
\begin{equation}
\label{eq:3d-Dx}
  \boxed{%
    \frac{\partial \bm{u}}{\partial x}
    = (\Iz \kron \Iy \kron \Dx)\,\bm{u}.
  }
\end{equation}
Reading this equation as an algorithm: for each pair $(j,k)$, apply
the small banded $\Dx$ to the $x$-line $u(:,j,k)$.  The outer identity
factors $\Iy$ and $\Iz$ say the $y$ and $z$ indices are passengers
during an $x$-sweep; they come along for the ride.  The $y$- and
$z$-derivatives are obtained by sliding the derivative block into
the appropriate slot,
\begin{equation}
\label{eq:3d-Dy-Dz}
  \frac{\partial \bm{u}}{\partial y}
  = (\Iz \kron \Dy \kron \Ix)\,\bm{u},
  \qquad
  \frac{\partial \bm{u}}{\partial z}
  = (\Dz \kron \Iy \kron \Ix)\,\bm{u},
\end{equation}
and the rest of the pattern follows.  A full sweep in any direction
touches each of the $\Nall$ grid values a fixed number of times, so
the cost is the optimal $\bigO(\Nall)$.

\subsection{Second derivatives}
\label{sec:3d-second}

The second derivative needs no new ideas.  Swap $\Dx$ for the
tridiagonal $\Dxx$ and do it again:
\begin{equation}
\label{eq:3d-Dxx}
  \frac{\partial^2 \bm{u}}{\partial x^2}
  = (\Iz \kron \Iy \kron \Dxx)\,\bm{u},
  \quad
  \frac{\partial^2 \bm{u}}{\partial y^2}
  = (\Iz \kron \Dyy \kron \Ix)\,\bm{u},
  \quad
  \frac{\partial^2 \bm{u}}{\partial z^2}
  = (\Dzz \kron \Iy \kron \Ix)\,\bm{u}.
\end{equation}
The kernel is different.  The algorithm is the same.

\subsection{Mixed derivatives}
\label{sec:3d-mixed}

The mixed derivative $\partial^2 u / \partial x\,\partial y$ is
two sweeps in a row, one per direction.  The mixed-product rule
(\cref{app:mixed-product}) collapses the composition into a single
Kronecker product:
\begin{equation}
\label{eq:3d-Dxy}
  \frac{\partial^2 \bm{u}}{\partial x\,\partial y}
  =
  \underbrace{(\Iz \kron \Dy \kron \Ix)}_{\text{then } y}
  \underbrace{(\Iz \kron \Iy \kron \Dx)}_{\text{first } x}\,\bm{u}
  \;=\;
  (\Iz \kron \Dy \kron \Dx)\,\bm{u}.
\end{equation}
The right-hand side is cosmetic.  The code still executes two
sweeps, an intermediate field between them, and never builds the
dense product $\Dy \kron \Dx$.

\subsection{The Laplacian}
\label{sec:3d-laplacian}

The seven-point discrete Laplacian at the interior point $(i,j,k)$,
\begin{equation}
\label{eq:laplacian-pointwise}
  (\Lap\,u)_{ijk}
  \;=\;
  \frac{u_{i-1,j,k} {-} 2u_{ijk} {+} u_{i+1,j,k}}{h_x^2}
  + \frac{u_{i,j-1,k} {-} 2u_{ijk} {+} u_{i,j+1,k}}{h_y^2}
  + \frac{u_{i,j,k-1} {-} 2u_{ijk} {+} u_{i,j,k+1}}{h_z^2},
\end{equation}
is the sum of three single-direction second derivatives.  Each term
looks at only one coordinate.  That is why the Laplacian falls
immediately into three additive sweeps,
\begin{equation}
\label{eq:3d-laplacian}
    \Lap\,\bm{u}
    \;=\;
    (\Iz \kron \Iy \kron \Dxx)\,\bm{u}
    \,+\,
    (\Iz \kron \Dyy \kron \Ix)\,\bm{u}
    \,+\,
    (\Dzz \kron \Iy \kron \Ix)\,\bm{u},
\end{equation}
which is exactly the identity that powered the introduction.  The
three results are computed line by line and added at the end.

\subsection{Memory layout and sweep direction}
\label{sec:memory}

An $x$-sweep is effortless because consecutive $x$-points sit next
to each other in memory.  A $y$-sweep strides across $N_x$ elements
to reach the next grid point in its line, and a $z$-sweep strides
across $N_x N_y$.  Two remedies exist: walk the stride with gather
and scatter, or transpose the active direction into the fast axis
before sweeping and transpose back afterward.  On modern caches and
accelerators the transpose almost always wins.  That is the serial
picture.  Its distributed counterpart is the pencil decomposition
(\cref{sec:pencil}), in which the transpose becomes a collective
MPI operation between global orientations of the field.  In both
cases the 1D kernel is unchanged; only the data motion around it
moves.

\section{Compact (Pad\'{e}) schemes}
\label{sec:compact}

\begin{figure}[t]
  \centering
  \includegraphics[width=0.62\linewidth]{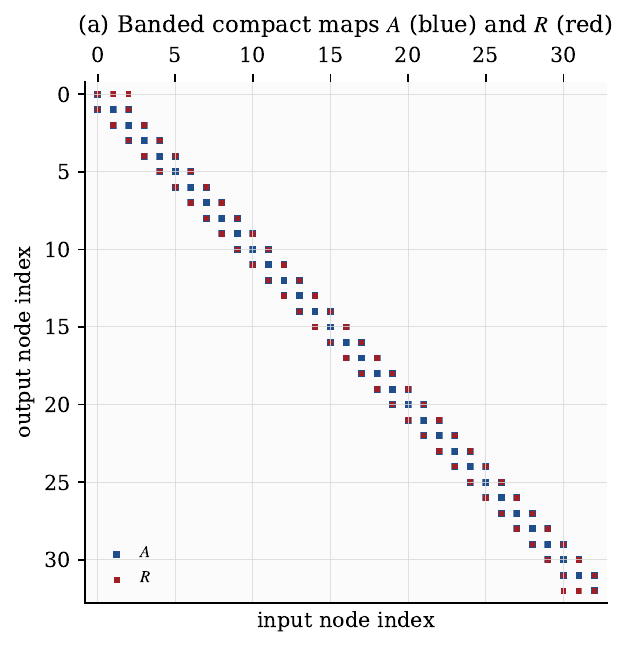}
  \caption{Compact first-derivative maps in factorized form.  The
  left matrix $\Ax$ and the right matrix $\Rx$ are both banded.  This
  is the important object for computation.  The operative algorithm is
  ``banded right-hand side plus banded line solve,'' not ``dense
  implicit differentiation matrix.''}
  \label{fig:compact-sparsity}
\end{figure}

Compact schemes trade a wider explicit stencil for a narrower implicit
relation: instead of writing the derivative as a long sum of neighbor
values, they couple a few nearby derivative values
together~\citep{lele1992compact,beam1976implicit}.  The reward is
spectral-like accuracy with a three-point stencil.  The price, or what
looks like a price on paper, is that one can no longer write
$\bm{u}' = \Dx\, \bm{u}$ with a banded $\Dx$.  Instead, one writes
\begin{equation*}
  \Ax\, \bm{u}' \;=\; \Rx\, \bm{u},
\end{equation*}
with $\Ax$ and $\Rx$ both banded.  The algebraically equivalent
statement $\bm{u}' = \Ax^{-1}\Rx\, \bm{u}$ conceals a pitfall:
$\Ax^{-1}\Rx$ is generically dense.  \Cref{fig:compact-sparsity}
shows why the factorized form is the only one that should ever touch
memory.  Both matrices are tridiagonal.  Keep them that way: apply
$\Rx$ (banded matvec, $\bigO(N_x)$), then solve with $\Ax$
(Thomas algorithm~\citep{thomas1949elliptic}, $\bigO(N_x)$).  Total
work per line: $\bigO(N_x)$.  Total work per 3D sweep: $\bigO(\Nall)$.
The compact scheme has the same asymptotic cost as central
differences, with a larger but controlled constant and typically much
better resolution per grid point.

\paragraph{Example: fourth-order Pad\'{e} first derivative.}
The classical fourth-order Pad\'{e} first derivative in 1D
reads~\citep{lele1992compact}
\begin{equation}
\label{eq:pade-example}
  \tfrac{1}{4}\,u'_{i-1} + u'_i + \tfrac{1}{4}\,u'_{i+1}
  =
  \tfrac{3}{2}\,\frac{u_{i+1} - u_{i-1}}{2h}.
\end{equation}
The two coefficients $\alpha = 1/4$ and $a = 3/2$ are determined by
matching the Taylor expansion through fourth order; the same
construction yields the sixth-order tridiagonal scheme with
$\alpha = 1/3$, $a = 14/9$, and $b = 1/9$ on a wider right-hand-side
stencil~\citep{lele1992compact}.  The left-hand side couples three
derivative values, hence the term ``compact,'' while the right-hand
side uses the standard centered stencil with a modified coefficient.
In matrix form this is $\Ax\,\bm{u}' = \Rx\,\bm{u}$, where $\Ax$ is
tridiagonal with $1$ on the main diagonal and $1/4$ on the sub- and
super-diagonals, and $\Rx$ encodes the right-hand side stencil.  The
equivalent expression $\bm{u}' = \Ax^{-1}\Rx\,\bm{u}$ is algebraically
correct, but the implementation proceeds by applying $\Rx$ and
solving with $\Ax$, rather than by forming $\Ax^{-1}$.

\paragraph{Three-dimensional form.}
In three dimensions, the compact scheme for $\partial u/\partial x$ at
each point $(i,j,k)$ reads
\begin{equation}
\label{eq:compact-3d-pointwise}
  \tfrac{1}{4}\,u'_{i-1,j,k} + u'_{ijk} + \tfrac{1}{4}\,u'_{i+1,j,k}
  =
  \tfrac{3}{2}\,\frac{u_{i+1,j,k} - u_{i-1,j,k}}{2h_x}.
\end{equation}
Written for all $N_x N_y N_z$ points simultaneously, this is the
system $\mathcal{A}\,\bm{u}' = \mathcal{R}\,\bm{u}$, where
$\mathcal{A}$ and $\mathcal{R}$ are $\Nall \times \Nall$ matrices.
Solving this as a monolithic system would be prohibitively expensive.
But because the coupling on both sides is only along $x$ (the indices
$j$ and $k$ are the same on both sides of the equation), the system
decomposes.

\paragraph{Kronecker form of compact schemes.}
The 3D compact first derivative in $x$ factors as
\begin{equation}
\label{eq:3d-compact-Dx}
  (\Iz \kron \Iy \kron \Ax)\;\bm{v}
  =
  (\Iz \kron \Iy \kron \Rx)\;\bm{u},
\end{equation}
where $\bm{v}$ is the vector of derivative values.  Because $\Ax$
appears inside the same Kronecker product structure, this system
decomposes into $N_y \cdot N_z$ independent tridiagonal solves of
size $N_x$, one per $x$-line.  Each line costs $\bigO(N_x)$, so the
full sweep costs $\bigO(\Nall)$.  The higher dimension does not
change the kernel.  It only increases the number of identical lines to
be processed.
Appendix~\ref{app:factorized-lines} writes this decomposition out
explicitly in $\vect{}$ form.

For a fixed $(j,k)$, the line operation is:
\begin{enumerate}
  \item Compute the right-hand side: $\bm{b} = \Rx\,u(:,j,k)$
        (banded matrix-vector product, $\bigO(N_x)$).
  \item Solve the tridiagonal system: $\Ax\,v(:,j,k) = \bm{b}$
        via the Thomas algorithm ($\bigO(N_x)$).
\end{enumerate}
Repeat for all $(j,k)$ pairs.  The $y$- and $z$-compact derivatives
are the same, with the appropriate direction:
\begin{equation}
  (\Iz \kron \Ay \kron \Ix)\;\bm{v}
  = (\Iz \kron \Ry \kron \Ix)\;\bm{u},
  \qquad
  (\Az \kron \Iy \kron \Ix)\;\bm{v}
  = (\Rz \kron \Iy \kron \Ix)\;\bm{u}.
\end{equation}
Compact second derivatives ($\Axx$, $\Rxx$, etc.) work identically.
The same factorized viewpoint applies.  The effective map
$\Axx^{-1}\Rxx$ may be dense when written explicitly, but the banded
apply remains linear-time because that dense product is never formed.

\paragraph{Non-uniform grids.}
Compact schemes can be derived on non-uniform grids by matching Taylor
expansions with the local spacing.  The matrices $\Ax$ and $\Rx$ get
point-dependent entries, but they remain tridiagonal, and the Kronecker
structure is unchanged.  The same factorized line solve survives on
stretched grids without any asymptotic penalty.

\section{Galerkin and B-spline methods}
\label{sec:galerkin}

Finite differences are not special.  Every Cartesian discretization
that uses a tensor-product basis factors the same way: a
finite-element method, a spectral-element method, a Fourier or
Chebyshev expansion, a B-spline collocation, an isogeometric
analysis.  The ingredients have different names, but the Kronecker
machinery is identical.  This section tracks the factorization
through the weak form.

\subsection{Galerkin on a tensor-product domain}
\label{sec:galerkin-general}

\begin{figure}[t]
  \centering
  \includegraphics[width=\linewidth]{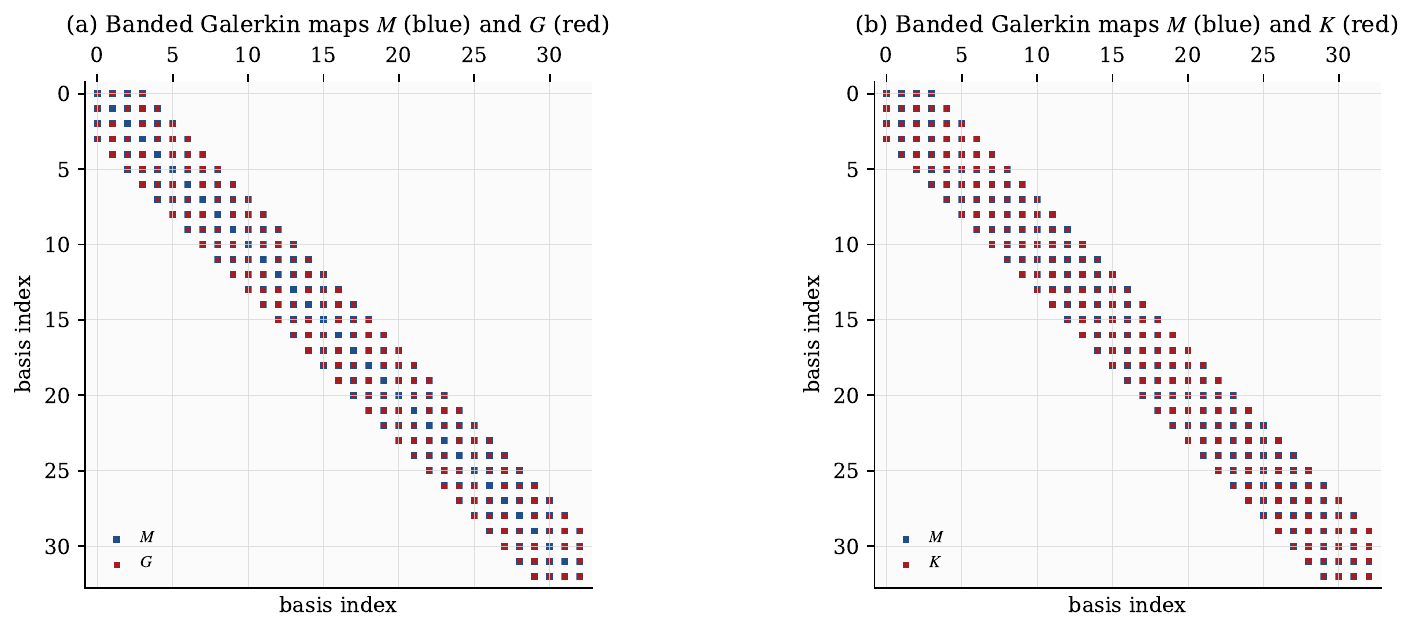}
  \caption{One-dimensional Galerkin spline maps.  The mass matrix $M$
  and the transport and stiffness matrices $G$ and $K$ are banded
  because each basis function overlaps only a fixed number of
  neighbors.  This is the locality that the tensor-product algorithm
  exploits in 3D.}
  \label{fig:galerkin-sparsity}
\end{figure}

On the box $\Omega = [a,b] \times [c,d] \times [e,f]$, the Galerkin
recipe for $-\nabla^2 u = g$ expands the solution in basis functions
$\{\Phi_n\}$, tests against the same family, and integrates by parts
to produce $K \hat{\bm{u}} = \bm{g}$ with
\begin{equation}
\label{eq:3d-mass-full}
  M_{mn} = \int_\Omega \Phi_m \Phi_n \, d\bm{x},
  \qquad
  K_{mn} = \int_\Omega \nabla\Phi_m \cdot \nabla\Phi_n \, d\bm{x}.
\end{equation}
Written at face value, $M$ and $K$ are $\Nall \times \Nall$.
Written with eyes open, they are not.  The 1D Galerkin building
blocks are banded (\cref{fig:galerkin-sparsity}), and the 3D
operators inherit that locality through a Kronecker factorization
that is about to fall out of the geometry of the basis.

\paragraph{The 3D operators factor.}
On a tensor-product domain with a tensor-product basis, every 3D
basis function splits as
\begin{equation}
\label{eq:galerkin-basis}
  \Phi_{ijk}(x,y,z) = \phi^x_i(x)\,\phi^y_j(y)\,\phi^z_k(z),
\end{equation}
and $d\bm{x}$ factors with it.  The 3D mass integral is therefore
the product of three 1D integrals,
\begin{equation}
\label{eq:mass-factor}
  M_{(ijk),(pqr)}
  = \underbrace{\int \phi^x_i\,\phi^x_p\,dx}_{(M_x)_{ip}}
  \;\underbrace{\int \phi^y_j\,\phi^y_q\,dy}_{(M_y)_{jq}}
  \;\underbrace{\int \phi^z_k\,\phi^z_r\,dz}_{(M_z)_{kr}}.
\end{equation}
This is exactly the Kronecker product:
\begin{equation}
\label{eq:3d-mass}
  \boxed{%
    M = M_z \kron M_y \kron M_x
  }
\end{equation}

For the stiffness matrix, the gradient introduces a derivative in one
direction at a time.  Taking the $x$-contribution as an example:
\begin{equation}
  \int_\Omega \frac{\partial\Phi_{ijk}}{\partial x}
  \frac{\partial\Phi_{pqr}}{\partial x}\,d\bm{x}
  = \underbrace{\int \phi'^x_i\,\phi'^x_p\,dx}_{(K_x)_{ip}}
  \;\underbrace{\int \phi^y_j\,\phi^y_q\,dy}_{(M_y)_{jq}}
  \;\underbrace{\int \phi^z_k\,\phi^z_r\,dz}_{(M_z)_{kr}}.
\end{equation}
Summing over all three gradient components gives the full stiffness
matrix:
\begin{equation}
\label{eq:3d-stiffness}
  K = M_z \kron M_y \kron K_x
    + M_z \kron K_y \kron M_x
    + K_z \kron M_y \kron M_x.
\end{equation}
The Galerkin semi-discrete system for the time-dependent problem
$\partial u/\partial t = \nabla^2 u$ is
$M\,\dot{\bm{u}} = -K\,\bm{u} + \bm{f}$, and it inherits full
Kronecker structure.

\paragraph{Mass inversion and stiffness apply.}
An explicit step needs $M^{-1}$ at every stage.  The inverse
rule~\eqref{eq:kron-inv} gives it without a fight:
\begin{equation}
  M^{-1} = M_z^{-1} \kron M_y^{-1} \kron M_x^{-1},
\end{equation}
and applying $M^{-1}$ is three passes of 1D banded solves along the
three line families.  Each $M_\xi$ is tridiagonal for linear
elements, pentadiagonal for quadratics, and banded with bandwidth
$2p+1$ for splines of degree $p$, so every pass is $\bigO(N_\xi)$
per line and $\bigO(\Nall)$ across the grid.  The stiffness apply
$K\bm{u}$ is the same story played three times, once for each term
in \cref{eq:3d-stiffness}; never form the Kronecker product
explicitly, apply each factor along its own direction and add.

\paragraph{Spectral methods are the same story.}
With global polynomials (Chebyshev, Legendre) or Fourier modes, the
1D mass matrices become diagonal in an orthogonal basis and the
stiffness matrices keep their Kronecker-sum shape.  Spectral element
methods restrict the same construction to each tensor-product
element and reap the same line
structure~\citep{deville2002high,haidvogel1979accurate}.  None of
the algebra changes; only the 1D pieces become polynomial rather
than finite-difference.

\subsection{B-spline and isogeometric discretizations}
\label{sec:bspline}

B-splines are the discretization the Kronecker picture was waiting
for~\citep{deboor1978practical}.  In 1D, a spline of degree $p$ is a
linear combination $u(x) = \sum_i B^x_i(x)\,\alpha_i$ of compactly
supported basis functions, and the coefficient vector $\alpha$ is
the object one stores.  Evaluating the spline and its derivatives is
a linear map:
\begin{equation}
  u(x) = B(x)\,\alpha,
  \qquad
  u_x(x) = B'(x)\,\alpha,
  \qquad
  u_{xx}(x) = B''(x)\,\alpha,
\end{equation}
the spline analog of a 1D differentiation matrix.  The evaluation
maps $B$, $B'$, and $B''$ have only $\bigO(p)$ nonzeros per row:
at any point, at most $p+1$ B-splines are active.  The Galerkin mass
and stiffness matrices assembled from these maps have bandwidth
$\bigO(p)$, commonly written as $2p+1$ for open knot vectors.  The
same banded structure appears in the boundary-consistent filtering
construction used for high-fidelity
turbulence~\citep{bay2019dissertation}, which again replaces a
three-dimensional filter by repeated one-dimensional spline maps.  The
same tensor-product spline structure also governs the statistical
problem: used as a regression basis rather than a discretization, it
makes the optimal resolution solvable in closed form---its cost set by
interaction order rather than ambient dimension---instead of found by
search~\citep{bay2026spline}.

\begin{figure}[t]
  \centering
  \includegraphics[width=\linewidth]{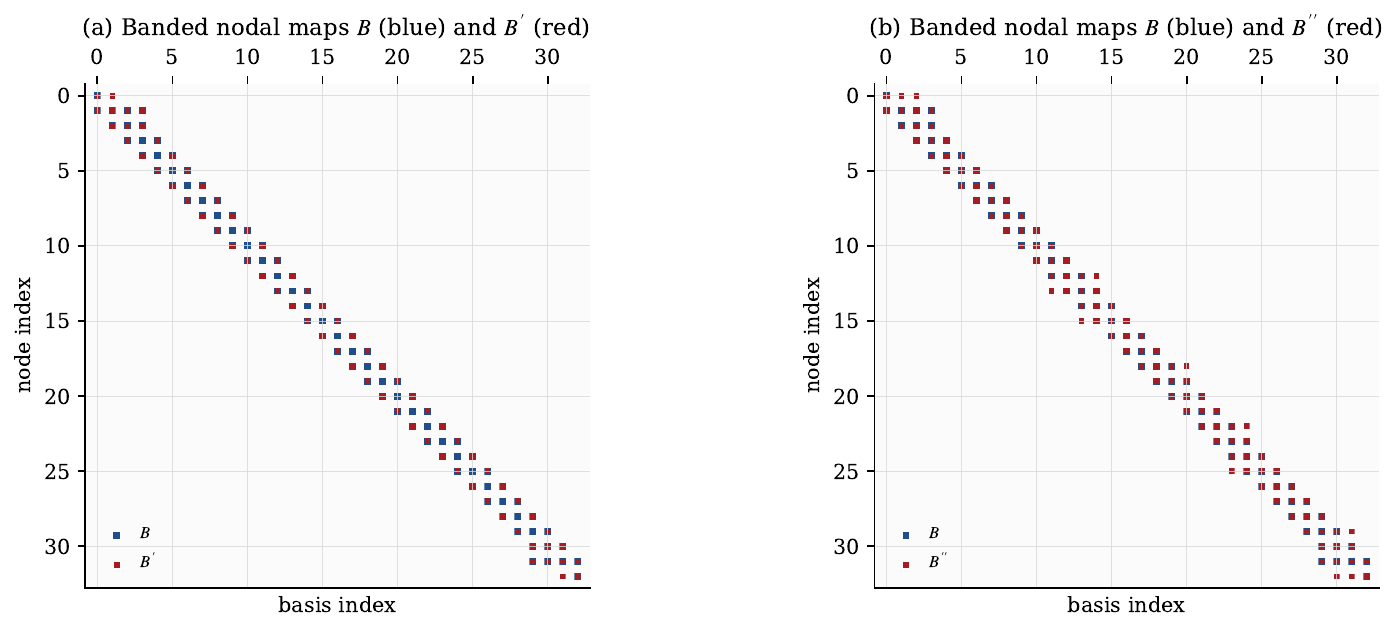}
  \caption{One-dimensional B-spline collocation maps.  The nodal maps
  $B$, $B'$, and $B''$ remain sparse and banded because compact
  support limits which basis functions can contribute at a point.
  Locality is carried by these one-dimensional maps.  The formed
  products $B'B^{-1}$ and $B''B^{-1}$ need not stay banded, so the
  efficient algorithm keeps the factorized form.}
  \label{fig:bspline-sparsity}
\end{figure}

Choose evaluation or quadrature points $\{x_q\}_{q=1}^{Q_x}$ and
form the matrices
\[
  (B_x)_{qi} = B^x_i(x_q),
  \qquad
  (B_x')_{qi} = \frac{dB^x_i}{dx}(x_q),
  \qquad
  (B_x'')_{qi} = \frac{d^2 B^x_i}{dx^2}(x_q).
\]
If $u \in \R^{Q_x}$ denotes the vector of sampled spline values
$u_q = u(x_q)$, then
\begin{equation}
  u = B_x \alpha,
  \qquad
  u_x = B_x' \alpha,
  \qquad
  u_{xx} = B_x'' \alpha.
\end{equation}
For Galerkin assembly, let $W_x = \diag(w_1,\ldots,w_{Q_x})$ be the
diagonal matrix of quadrature weights.  The 1D Galerkin matrices are
then
\begin{equation}
  M_x = B_x^T W_x B_x,
  \qquad
  G_x = B_x^T W_x B_x',
  \qquad
  K_x = (B_x')^T W_x B_x'.
\end{equation}
Here $M_x$ is the 1D mass matrix, $G_x$ is the first-derivative
coupling matrix, and $K_x$ is the 1D stiffness matrix.  The same
construction gives $M_y$, $M_z$, $G_y$, $G_z$, $K_y$, and $K_z$.  In
three dimensions,
\begin{equation}
  u(x,y,z)
  = \sum_{i=1}^{N_x}\sum_{j=1}^{N_y}\sum_{k=1}^{N_z}
  \alpha_{ijk}\, B^x_i(x)\, B^y_j(y)\, B^z_k(z),
\end{equation}
and the assembled operators are still
\begin{equation}
  M = M_z \kron M_y \kron M_x,
  \qquad
  K = M_z \kron M_y \kron K_x
    + M_z \kron K_y \kron M_x
    + K_z \kron M_y \kron M_x.
\end{equation}
The 3D spline solve therefore has the same tensor-product structure as
the finite-element case above.  The only additional ingredient is the
assembly of the 1D matrices.  Because each spline overlaps only
$\bigO(p)$ neighbors, the 1D mass and stiffness matrices remain banded
with bandwidth $2p+1$, and every line operation remains
$\bigO(N_\xi)$.  Consequently, the 3D Galerkin or isogeometric apply is
linear in $\Nall$ because it is built entirely from these banded 1D
blocks.

\paragraph{Isogeometric analysis.}
In isogeometric analysis
(IGA)~\citep{hughes2005isogeometric,cottrell2009isogeometric}, the
same B-spline (or NURBS) basis that represents the CAD geometry is
reused as the solution basis.  This gives exact CAD geometry on each
tensor-product patch, while preserving the same operator algebra.  On
a box-shaped parameter domain, the primary unknown is usually the
coefficient vector $\alpha$, not sampled values $u$.  If one wants the
\emph{values} of the derivative at sample points, the apply is direct:
\[
  u_x = B_x' \alpha.
\]
If instead one wants a derivative operator that maps spline
coefficients to spline coefficients in the same basis, introduce
coefficients $\beta_x$ such that $B(x)\,\beta_x$ is the $L^2$
projection of $u_x$ back into the spline space.  The 1D projection
equation is
\begin{equation}
  M_x \beta_x = G_x \alpha,
  \qquad\text{so}\qquad
  \beta_x = D_x^{\mathrm{IGA}} \alpha,
  \qquad
  D_x^{\mathrm{IGA}} = M_x^{-1} G_x.
\end{equation}
This is the IGA analog of applying a 1D derivative matrix along a
line.  But the factorized form is the computational one that matters:
apply $G_x$, then solve with $M_x$.  If one forms $M_x^{-1}G_x$
explicitly, that product is generally dense.  In 3D, the same
operator factors as
\[
  I_z \kron I_y \kron D_x^{\mathrm{IGA}},
\]
and similarly in $y$ and $z$.  For diffusion and Poisson problems, the
second-derivative action usually enters through the stiffness matrix
$K_x$ rather than an explicit pointwise $\Dxx$.

The computational consequence is therefore straightforward: IGA changes
the basis and often improves accuracy per degree of freedom through
higher continuity, but it does not alter the 3D apply.  Once the 1D
spline matrices are assembled, every matrix-vector product or implicit
solve remains a batch of 1D line operations.

\paragraph{B-spline collocation.}
An alternative to the Galerkin approach is B-spline collocation:
instead of integrating against test functions, evaluate the PDE
directly at a set of collocation points and require the B-spline
expansion to satisfy the equation pointwise.  A common choice is the
Greville points~\citep{johnson2005higher}, but the tensor algebra does
not depend on the specific set of collocation sites.  In 1D, evaluate
the basis and its derivatives at the collocation points to obtain
square matrices $B_x$, $B_x'$, and $B_x''$.  If $u$ denotes the vector
of spline values at those points, then
\begin{equation}
  u = B_x \alpha,
  \qquad
  u_x = B_x' \alpha,
  \qquad
  u_{xx} = B_x'' \alpha.
\end{equation}
Eliminating the coefficients gives the collocation differentiation
matrices
\begin{equation}
  \Dx = B_x' B_x^{-1},
  \qquad
  \Dxx = B_x'' B_x^{-1}.
\end{equation}
These are the spline-induced first- and second-derivative matrices on
collocated values.  They are useful for analysis, but they are not the
best objects to form in code.  The matrices $B_x$, $B_x'$, and
$B_x''$ in \Cref{fig:bspline-sparsity} are banded, whereas
$B_x^{-1}$ is generally dense.  So the linear-time collocation apply
should be read in factorized form: solve $B_x \alpha = u$ on each
line, then apply $B_x'$ or $B_x''$ to that recovered coefficient line.
Appendix~\ref{app:factorized-lines} writes the same coefficient-recovery
step directly in tensor-product form.
This banded solve plus banded matvec still costs $\bigO(N_x)$ per line
and therefore $\bigO(\Nall)$ over the full 3D grid.
In periodic shift-invariant settings, these formed matrices can become
circulant and enjoy useful symmetry properties, but they are still not
generically banded once $B_x^{-1}$ is formed.

In coefficient form, the 3D collocation Poisson operator can be
written without ever introducing a dense value-space differentiation
matrix:
\begin{equation}
  -\bigl(
    B_z \kron B_y \kron B_x''
    + B_z \kron B_y'' \kron B_x
    + B_z'' \kron B_y \kron B_x
  \bigr)\,\alpha
  = \bm{g}.
\end{equation}
This representation is the cleanest implementation form because every
factor on every line is banded.

For the Poisson equation $-\nabla^2 u = g$, the full 3D collocation
system on collocated values can also be written as
\begin{equation}
  -\bigl(\Iz \kron \Iy \kron \Dxx
       + \Iz \kron \Dyy \kron \Ix
       + \Dzz \kron \Iy \kron \Ix\bigr)\,\bm{u}
  = \bm{g},
\end{equation}
where $\Dxx$, $\Dyy$, $\Dzz$ are the spline collocation
second-derivative matrices in each direction.  This is a Kronecker
sum, identical in structure to the finite-difference Laplacian
\eqref{eq:3d-laplacian}, and it can be solved by the same fast
diagonalization technique (\cref{sec:fast-diag}).

Collocation avoids numerical quadrature entirely, which simplifies
setup.  The computational message is the same as for compact
schemes: use the factorized banded maps, not a preformed dense
derivative matrix, and the per-line operations remain exactly the same
\texttt{sweep} and \texttt{solve} from the pseudocode
(\cref{fig:pseudocode}).

\section{Implicit time stepping via ADI}
\label{sec:implicit}

Stiffness forces an implicit treatment.  For the heat equation on a
fine grid, the explicit stability limit $\dt \lesssim h^2/\nu$ is
prohibitive, and implicit methods lift that ceiling by solving a
system at each step.  In three dimensions that system is
$\Nall \times \Nall$.  Forming it is out of the question.  The
classical cure, invented in the 1950s for oil-reservoir simulation by
\citet{peaceman1955numerical} and generalized to three space variables
by \citet{douglas1962alternating,douglas1964numerical}, is the
alternating-direction implicit (ADI) splitting: replace one coupled 3D
solve by three 1D solves, one per direction, each of which is a batch
of banded line solves.

\begin{figure}[t]
  \centering
  \resizebox{\linewidth}{!}{%
  \begin{tikzpicture}[
    x={(0.72cm,0cm)}, y={(0cm,0.72cm)}, z={(-0.48cm,-0.36cm)},
    every node/.style={font=\small},
  ]
    \pgfmathsetmacro{\Nx}{4}
    \pgfmathsetmacro{\Ny}{3}
    \pgfmathsetmacro{\Nz}{3}
    \definecolor{sweepblue}{HTML}{1f4e8a}
    \definecolor{sweepred}{HTML}{a11e24}
    \definecolor{sweepgreen}{HTML}{1f7a3a}
    \definecolor{boxedge}{HTML}{7a7a7a}
    \definecolor{faintline}{HTML}{c8c8c8}

    \newcommand{\wireframe}{%
      \draw[boxedge, line width=0.5pt] (0,0,0) -- (\Nx,0,0) -- (\Nx,\Ny,0) -- (0,\Ny,0) -- cycle;
      \draw[boxedge, line width=0.5pt] (0,0,\Nz) -- (\Nx,0,\Nz) -- (\Nx,\Ny,\Nz) -- (0,\Ny,\Nz) -- cycle;
      \draw[boxedge, line width=0.5pt] (0,0,0) -- (0,0,\Nz);
      \draw[boxedge, line width=0.5pt] (\Nx,0,0) -- (\Nx,0,\Nz);
      \draw[boxedge, line width=0.5pt] (\Nx,\Ny,0) -- (\Nx,\Ny,\Nz);
      \draw[boxedge, line width=0.5pt] (0,\Ny,0) -- (0,\Ny,\Nz);
    }

    \begin{scope}[xshift=0cm]
      \node[anchor=south, font=\bfseries, sweepblue]
        at (\Nx/2,\Ny+2.9,\Nz) {Step 1.\ solve along $x$};
      \node[anchor=south] at (\Nx/2,\Ny+2.0,\Nz)
        {$(I-\tfrac{\nu\dt}{2}L_x)\,\bm{u}^{*}=\ldots$};
      \wireframe
      \foreach \k in {0,...,\Nz}{
        \foreach \j in {0,...,\Ny}{
          \draw[faintline, line width=0.45pt] (0,\j,\k) -- (\Nx,\j,\k);
        }
      }
      \draw[sweepblue, line width=1.9pt, -{Stealth[length=2.8mm]}]
        (-0.25,1,1) -- (\Nx+0.25,1,1);
    \end{scope}

    \begin{scope}[x={(1cm,0cm)}, y={(0cm,1cm)}]
      \draw[line width=1.4pt, -{Stealth[length=3.5mm, width=3.0mm]}, black!75]
        (3.25,1.0) -- (4.25,1.0);
    \end{scope}

    \begin{scope}[xshift=7.0cm]
      \node[anchor=south, font=\bfseries, sweepred]
        at (\Nx/2,\Ny+2.9,\Nz) {Step 2.\ solve along $y$};
      \node[anchor=south] at (\Nx/2,\Ny+2.0,\Nz)
        {$(I-\tfrac{\nu\dt}{2}L_y)\,\bm{u}^{**}=\ldots$};
      \wireframe
      \foreach \k in {0,...,\Nz}{
        \foreach \i in {0,...,\Nx}{
          \draw[faintline, line width=0.45pt] (\i,0,\k) -- (\i,\Ny,\k);
        }
      }
      \draw[sweepred, line width=1.9pt, -{Stealth[length=2.8mm]}]
        (2,-0.25,1) -- (2,\Ny+0.25,1);
    \end{scope}

    \begin{scope}[x={(1cm,0cm)}, y={(0cm,1cm)}]
      \draw[line width=1.4pt, -{Stealth[length=3.5mm, width=3.0mm]}, black!75]
        (10.25,1.0) -- (11.25,1.0);
    \end{scope}

    \begin{scope}[xshift=14.0cm]
      \node[anchor=south, font=\bfseries, sweepgreen]
        at (\Nx/2,\Ny+2.9,\Nz) {Step 3.\ solve along $z$};
      \node[anchor=south] at (\Nx/2,\Ny+2.0,\Nz)
        {$(I-\tfrac{\nu\dt}{2}L_z)\,\bm{u}^{n+1}=\ldots$};
      \wireframe
      \foreach \j in {0,...,\Ny}{
        \foreach \i in {0,...,\Nx}{
          \draw[faintline, line width=0.45pt] (\i,\j,0) -- (\i,\j,\Nz);
        }
      }
      \draw[sweepgreen, line width=1.9pt, -{Stealth[length=2.8mm]}]
        (2,1,-0.25) -- (2,1,\Nz+0.25);
    \end{scope}
  \end{tikzpicture}}
  \caption{The ADI cascade.  A single implicit step on the 3D heat
  equation is replaced by three sequential line solves: first along
  every $x$-line, then along every $y$-line, then along every $z$-line.
  Each substep is a batch of $\bigO(N_\xi)$ banded line solves, so the
  entire implicit step is $\bigO(\Nall)$, the same asymptotic cost as
  an explicit sweep, with the explicit stability restriction gone.}
  \label{fig:adi-cascade}
\end{figure}
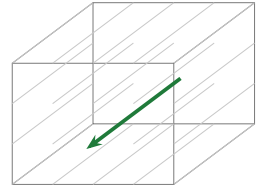

\begin{figure}[t]
  \centering
  \includegraphics[width=\linewidth]{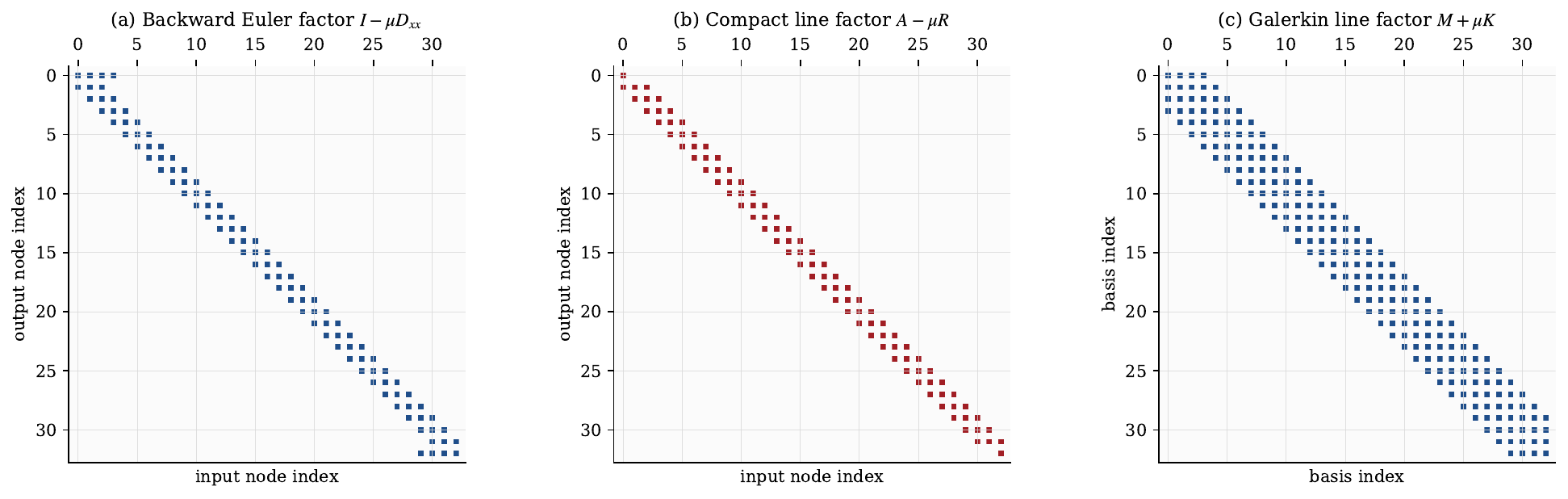}
  \caption{Typical one-dimensional line factors for implicit time
  stepping.  Whether the underlying discretization is explicit
  finite difference, compact Pad\'{e}, or Galerkin, the ADI substep is
  carried by a banded 1D solve.}
  \label{fig:implicit-sparsity}
\end{figure}

\Cref{fig:adi-cascade} shows the cascade pictorially: the coupled 3D
solve unravels, direction by direction, into banded line solves.
\Cref{fig:implicit-sparsity} shows the one-dimensional factors
themselves: explicit finite differences give backward Euler factors
$I - \tfrac{\nu\dt}{2}\Dxx$ directly, compact schemes give banded
factors $\Axx - \tfrac{\nu\dt}{2}\Rxx$ after moving the compact mass
to the operator side, and Galerkin and spline methods give factors
$M_x + \tfrac{\nu\dt}{2}K_x$ that are banded because their 1D
ingredients are.

\subsection{The splitting}

For the heat equation
$\partial u/\partial t = \nu \nabla^2 u$ with diffusivity $\nu$,
the semi-discrete system is $\dot{\bm{u}} = \nu\Lap\bm{u}$ and a
backward-Euler step is
\begin{equation}
\label{eq:backward-euler}
  (I - \nu\,\dt\,\Lap)\;\bm{u}^{n+1} = \bm{u}^n.
\end{equation}
The matrix on the left is $\Nall \times \Nall$.  Factoring it head-on
would defeat the whole enterprise.  The trick is to remember that
$\Lap = L_x + L_y + L_z$ with
\begin{equation}
  L_x = \Iz \kron \Iy \kron \Dxx, \quad
  L_y = \Iz \kron \Dyy \kron \Ix, \quad
  L_z = \Dzz \kron \Iy \kron \Ix,
\end{equation}
and to solve with each piece in turn.  The Douglas--Gunn
splitting~\citep{douglas1964numerical}, written here in its
Crank--Nicolson form, is the scheme of choice.  Set $r = \nu\dt/2$;
the three substeps are
\begin{align}
  (I - r\,L_x)\;\bm{u}^*
  &= \bigl(I + r\,L_x + 2r\,L_y + 2r\,L_z\bigr)\,\bm{u}^n,
  \label{eq:adi-x} \\
  (I - r\,L_y)\;\bm{u}^{**}
  &= \bm{u}^* - r\,L_y\,\bm{u}^n,
  \label{eq:adi-y} \\
  (I - r\,L_z)\;\bm{u}^{n+1}
  &= \bm{u}^{**} - r\,L_z\,\bm{u}^n.
  \label{eq:adi-z}
\end{align}
The first-step right-hand side $(I + r L_x + 2r(L_y + L_z))$ is the
Douglas--Gunn correction that makes the three substeps combine to
the full Crank--Nicolson update $(I - r\Lap)\,\bm{u}^{n+1} = (I + r\Lap)\,\bm{u}^n$
modulo a splitting error of order $\bigO((\nu\dt)^3)$ per
step~\citep{douglas1964numerical}; dropping the $2r(L_y + L_z)\,\bm{u}^n$
correction collapses the scheme to $\bigO(\dt)$ accuracy and biases
the heat operator toward the $x$-direction.  Each substep is a banded
1D solve because, under the identity from \cref{app:vec} plus
distributivity, the $x$-substep operator is
\begin{equation}
  I - \tfrac{\nu\dt}{2}\,L_x
  \;=\;
  \Iz \kron \Iy \kron \bigl(\Ix - \tfrac{\nu\dt}{2}\,\Dxx\bigr),
\end{equation}
and the bracketed 1D factor is tridiagonal.  The $y$- and
$z$-substeps factor identically in their own directions.  Each
substep therefore costs $\bigO(\Nall)$ in the Thomas algorithm, so
the whole implicit step is $\bigO(\Nall)$ with a modest constant,
the same asymptotic cost as one explicit sweep.

\subsection{The same splitting for every discretization}

Changing the discretization changes the $1$D factor, not the outer
Kronecker structure.  A sixth-order compact scheme, where
$\Dxx = \Axx^{-1}\Rxx$, produces the line factor
\[
  \Axx - \tfrac{\nu\dt}{2}\Rxx,
\]
tridiagonal because tridiagonals subtract.  A Galerkin or B-spline
semi-discretization $M\dot{\bm{u}} = -K\bm{u}$ produces
\[
  M_x + \tfrac{\nu\dt}{2} K_x,
\]
banded because both $M_x$ and $K_x$ are banded.  B-spline
collocation in coefficient form produces
\[
  B_x - \tfrac{\nu\dt}{2}\,B_x'',
\]
banded because both $B_x$ and $B_x''$ are.  In every case each ADI
substep is a banded 1D solve; only the bandwidth and the entries
change.  The same splitting also covers the
advection-diffusion equation
$\partial_t u + \bm{c}\cdot\nabla u = \nu\nabla^2 u$
with $L_\xi = c_\xi D_\xi + \nu D_{\xi\xi}$, and any equation whose
spatial operator decomposes cleanly along the coordinate axes.

\section{Direct Poisson solvers via fast diagonalization}
\label{sec:fast-diag}

There is a second, complementary payoff of the Kronecker structure.
For constant-coefficient separable problems, in particular the Poisson
equation $\Lap\,\bm{u} = \bm{f}$ and shifted Helmholtz equations of
the form $(\alpha I + \beta\Lap)\,\bm{u} = \bm{f}$ that arise inside
implicit time stepping, the same algebra gives a \emph{direct} solver
with no iteration at all.  The idea is due to
\citet{lynch1964direct}: if the one-dimensional second-derivative
matrices can be diagonalized, then diagonalizing each direction in
turn decouples the 3D problem into independent scalar divisions.  For
Chebyshev expansions, the companion construction is
\citet{haidvogel1979accurate}.  The eigenbases themselves are computed
once as a preprocessing step and reused on every solve, so the
expensive part happens only at setup time.

\subsection{Diagonalize each direction, add the eigenvalues}

Suppose each 1D second-derivative matrix is diagonalizable,
\begin{equation}
  \Dxx = S_x \Lambda_x S_x^{-1},
  \quad
  \Dyy = S_y \Lambda_y S_y^{-1},
  \quad
  \Dzz = S_z \Lambda_z S_z^{-1},
\end{equation}
with $\Lambda_\xi = \diag(\lambda_1^\xi, \ldots, \lambda_{N_\xi}^\xi)$.
Substituting into the Laplacian, applying the mixed-product rule
three times, and transforming $\bm{u}$ and $\bm{f}$ into the
tensor-product eigenbasis
$\hat{\bm{u}} = (S_z^{-1} \kron S_y^{-1} \kron S_x^{-1})\bm{u}$
leaves the Poisson problem diagonal
(\cref{app:fast-diag-derivation}):
\begin{equation}
\label{eq:poisson-diag}
  \boxed{%
    \hat{u}_{ijk}
    \;=\;
    \frac{\hat{f}_{ijk}}{\lambda_i^x + \lambda_j^y + \lambda_k^z}.
  }
\end{equation}
The 3D Poisson solve has become a scalar division at every grid
point.  Adding the three eigenvalues is the whole algorithm.  When a
denominator is zero, as in a pure Neumann Poisson problem, the usual
compatibility condition must hold and the zero mode is fixed by a
normalization such as zero mean.

\begin{figure}[t]
  \centering
  \resizebox{0.92\linewidth}{!}{%
  \begin{tikzpicture}[
    every node/.style={font=\small},
    box/.style={draw=black!65, rounded corners=2pt, minimum width=2.6cm, minimum height=0.9cm, align=center, fill=black!3},
    arr/.style={-{Stealth[length=3mm]}, line width=0.9pt, black!70}
  ]
    \node[box] (f) at (0,0) {$f$ on grid};
    \node[box] (tx) at (3.3,0) {$S_x^{-1}$ on\\$x$-lines};
    \node[box] (ty) at (6.6,0) {$S_y^{-1}$ on\\$y$-lines};
    \node[box] (tz) at (9.9,0) {$S_z^{-1}$ on\\$z$-lines};
    \node[box, fill=black!6] (div) at (13.4,0) {divide by\\$\lambda_i^x+\lambda_j^y+\lambda_k^z$};
    \node[box] (it) at (16.9,0) {inverse transforms\\by direction};
    \node[box] (u) at (20.0,0) {$u$ on grid};
    \draw[arr] (f) -- (tx);
    \draw[arr] (tx) -- (ty);
    \draw[arr] (ty) -- (tz);
    \draw[arr] (tz) -- (div);
    \draw[arr] (div) -- (it);
    \draw[arr] (it) -- (u);
    \node[anchor=north, font=\footnotesize] at (10,-0.75)
      {fast diagonalization is still a sequence of one-dimensional line operations};
  \end{tikzpicture}}
  \caption{Fast diagonalization pipeline.  The transform to the tensor-product
  eigenbasis is applied one direction at a time, then the solve is a
  scalar division, and the inverse transform returns to physical space.
  The direct solver is not a 3D factorization; it is a sequence of 1D
  transforms wrapped around a pointwise diagonal solve.}
  \label{fig:fastdiag-pipeline}
\end{figure}
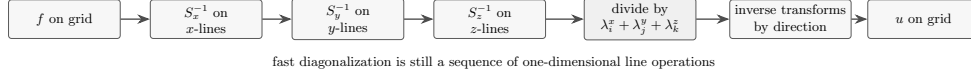

The bookkeeping around that scalar division is three transforms and
their inverses: apply $S_x^{-1}$ along all $x$-lines, then $S_y^{-1}$
along all $y$-lines, then $S_z^{-1}$ along all $z$-lines to get
$\hat{\bm{f}}$; divide pointwise; apply the forward transforms in
reverse (\cref{fig:fastdiag-pipeline}).  For a uniform grid the
eigenvectors are trigonometric and the boundary condition, grid
placement, and endpoint closure select the appropriate member of the
standard transform families: discrete sine transforms for homogeneous
Dirichlet boundaries, discrete cosine transforms for homogeneous
Neumann boundaries, and the discrete Fourier transform for periodic
boundaries.  Each of these costs $\bigO(N_\xi \log N_\xi)$ per line
and $\bigO(\Nall \log N_{\max})$ overall, and all are available in
optimized form (for instance through FFTW~\citep{frigo2005design}).  On stretched or Chebyshev
grids the eigenvectors are not generally trigonometric; the transforms
become small dense matrix products at $\bigO(N_\xi^2)$ per line and
$\bigO(\Nall N_{\max})$ overall, still far below the cost of factoring
the full 3D operator~\citep{haidvogel1979accurate}.  A shifted problem
$(\alpha I+\beta\Lap)\bm{u}=\bm{f}$ is handled at zero conceptual
extra cost by replacing the denominator in~\cref{eq:poisson-diag} with
$\alpha+\beta(\lambda_i^x+\lambda_j^y+\lambda_k^z)$, provided that
quantity never vanishes except for modes fixed by the compatibility
condition.

\section{Applicability of the tensor-product reduction}
\label{sec:applicability}

The right question is not whether a problem is ``three-dimensional.''
Every problem in this paper is three-dimensional.  The right question is
whether the algebra can be written as a small sum of tensor products,
\begin{equation}
\label{eq:applicability-kron-sum}
  L = \sum_{r=1}^{R}
  C_z^{(r)} \kron C_y^{(r)} \kron C_x^{(r)}.
\end{equation}
When $R$ is small and each one-dimensional factor is banded, diagonal,
or transformable, the whole machinery applies.  A single derivative
sweep has $R=1$ with two identity factors.  The Cartesian Laplacian has
$R=3$.  A separable coefficient such as
$a(x,y,z)=a_x(x)a_y(y)a_z(z)$ simply inserts diagonal one-dimensional
factors into the same product.

This criterion is the clean line between the exact reduction and the
many useful extensions of it.  Tensor-product grids and tensor-product
bases give the reduction natively.  Problems with additional
nonseparable physics still benefit from the same 1D solvers as split
operators, fast preconditioners, or patchwise kernels.  The table below
is therefore a map of where the reduction is exact and where it becomes
the computational backbone inside a larger solver.

\medskip
\footnotesize
\renewcommand{\arraystretch}{1.24}
\begin{tabular}{@{}>{\raggedright\arraybackslash}p{0.30\textwidth}>{\raggedright\arraybackslash}p{0.64\textwidth}@{}}
\toprule
\textbf{Setting} & \textbf{Role of the tensor-product reduction} \\
\midrule
\textbf{Uniform Cartesian} & \textbf{Exact.} The simplest case, with identical 1D
  factors reused on every line. \\
\addlinespace[5pt]
\textbf{Non-uniform Cartesian} & \textbf{Exact.} Spacing varies per direction, but the grid
  is still a tensor product of 1D grids. \\
\addlinespace[5pt]
\textbf{Stretched or clustered grids} & \textbf{Exact.} Wall clustering, tanh maps, and
  biased spacing only change the 1D weights. \\
\addlinespace[5pt]
\textbf{Mixed boundary conditions} & \textbf{Exact.} Different boundary closures are
  absorbed into the endpoint rows of the 1D operators. \\
\addlinespace[5pt]
\textbf{Galerkin tensor-product basis} & \textbf{Exact.} Mass and stiffness matrices inherit
  the Kronecker structure of the basis. \\
\addlinespace[5pt]
\textbf{B-spline and IGA patches} & \textbf{Exact on tensor-product patches.} The same
  factorization holds for spline mass, stiffness, and collocation maps. \\
\addlinespace[5pt]
\textbf{B-spline collocation} & \textbf{Exact in factorized form.} Use banded maps $B$,
  $B'$, and $B''$ rather than preformed dense products such as
  $B''B^{-1}$. \\
\addlinespace[5pt]
\textbf{Separable variable coefficients} & \textbf{Exact.} Products such as
  $a_x(x)a_y(y)a_z(z)$ preserve \eqref{eq:applicability-kron-sum}. \\
\addlinespace[5pt]
\textbf{General variable coefficients} & \textbf{Core solver component.} The separable or
  constant-coefficient part remains a fast split operator or
  preconditioner; the remaining coupling is handled matrix-free,
  iteratively, or by a low-rank separated approximation when such an
  approximation is accurate. \\
\addlinespace[5pt]
\textbf{Curvilinear grids} & \textbf{Exact when the metrics separate, otherwise patchwise.}
  General metric terms add couplings, while tensor-product blocks still
  supply fast local kernels and preconditioners. \\
\addlinespace[5pt]
\textbf{Unstructured grids} & \textbf{Indirect.} There are no global grid lines, but
  tensor-product elements, blocks, and patches retain the same
  dimension-by-dimension kernels locally. \\
\bottomrule
\end{tabular}
\normalsize
\medskip

\section{How production codes actually run: multi-RHS kernels, sum factorization, and pencils}
\label{sec:practice}

So far the paper has argued that a three-dimensional Cartesian operator
factors into a loop of one-dimensional line kernels.  That argument is
complete at the algebraic level.  But a graduate student writing a
research code should know the three practical tricks that separate a
textbook implementation of the same loop from a production one running
near peak floating-point rate.  These tricks are almost never written
down in one place.  They are the reason high-order CFD and
spectral-element codes achieve their reputations for efficiency.

\subsection{Trick 1: reshape a sweep into a multi-right-hand-side line kernel}
\label{sec:blas3}

Consider the $x$-sweep $v_{ijk} = \sum_{\ell} (\Dx)_{i\ell}\, u_{\ell j k}$.
Naively this is $N_y N_z$ separate line operations, each of size
$N_x$.  That is the right mathematical kernel but the wrong software
shape: the operator data, boundary rows, and factorizations are reused
across many lines, yet a scalar loop exposes only one right-hand side
at a time.

\begin{figure}[t]
  \centering
  \resizebox{0.95\linewidth}{!}{%
  \begin{tikzpicture}[every node/.style={font=\small}]
    \definecolor{blockblue}{HTML}{1f4e8a}
    \definecolor{blockred}{HTML}{a11e24}
    \definecolor{blockfill}{HTML}{eef2f7}
    \definecolor{bandfill}{HTML}{c9d6e6}

    \begin{scope}[shift={(0,0)}]
      \pgfmathsetmacro{\tW}{1.6}
      \pgfmathsetmacro{\tH}{1.6}
      \pgfmathsetmacro{\tD}{0.7}
      \draw[fill=blockfill!50, draw=black!50]
        (\tD,\tD) -- (\tD+\tW,\tD) -- (\tD+\tW,\tD+\tH) -- (\tD,\tD+\tH) -- cycle;
      \draw[fill=blockfill, draw=black!70]
        (0,0) -- (\tW,0) -- (\tW,\tH) -- (0,\tH) -- cycle;
      \draw[black!50] (0,0) -- (\tD,\tD);
      \draw[black!50] (\tW,0) -- (\tW+\tD,\tD);
      \draw[black!50] (\tW,\tH) -- (\tW+\tD,\tD+\tH);
      \draw[black!50] (0,\tH) -- (\tD,\tD+\tH);
      \foreach \y in {0.4,0.8,1.2}{
        \draw[blockblue!80, line width=0.4pt] (0,\y) -- (\tW,\y);
      }
      \foreach \x in {0.4,0.8,1.2}{
        \draw[blockblue!80, line width=0.4pt] (\x,0) -- (\x,\tH);
      }
      \node[anchor=north] at (\tW/2+\tD/2,-0.25) {tensor $u_{ijk}$};
      \node[anchor=north, font=\footnotesize]
        at (\tW/2+\tD/2,-0.85) {shape $(N_x,N_y,N_z)$};
    \end{scope}

    \draw[-{Stealth[length=3mm]}, line width=0.8pt]
      (2.7,0.9) -- (4.1,0.9);
    \node[anchor=south, font=\footnotesize] at (3.4,0.95) {reshape};

    \begin{scope}[shift={(4.5,0)}]
      \pgfmathsetmacro{\mW}{4.0}
      \pgfmathsetmacro{\mH}{1.6}
      \draw[fill=blockfill, draw=black!70] (0,0) rectangle (\mW,\mH);
      \foreach \y in {0.4,0.8,1.2}{
        \draw[blockblue!80, line width=0.4pt] (0,\y) -- (\mW,\y);
      }
      \foreach \x in {0.4,0.8,...,3.6}{
        \draw[blockblue!80, line width=0.4pt] (\x,0) -- (\x,\mH);
      }
      \node[anchor=north] at (\mW/2,-0.25) {matrix $U$};
      \node[anchor=north, font=\footnotesize]
        at (\mW/2,-0.85) {shape $(N_x,\ N_y N_z)$};
      \draw[blockred, line width=1.2pt] (1.2,0.05) -- (1.2,\mH-0.05);
      \node[blockred, font=\footnotesize, anchor=south]
        at (1.2,\mH+0.05) {one $x$-line};
    \end{scope}

    \node at (10.3,0.9) {$V \;=\; D_x \cdot U$};

    \begin{scope}[shift={(11.7,0)}]
      \pgfmathsetmacro{\dS}{1.6}   
      \pgfmathsetmacro{\nc}{8}     
      \pgfmathsetmacro{\dC}{\dS/\nc}
      \draw[fill=blockfill, draw=black!70] (0,0) rectangle (\dS,\dS);
      \foreach \i in {0,...,7}{
        \pgfmathsetmacro{\row}{\i}
        \foreach \j in {\i, -1+\i, 1+\i}{
          \pgfmathsetmacro{\jj}{int(\j)}
          \ifnum\jj>-1\ifnum\jj<8
            \pgfmathsetmacro{\xx}{\jj*\dC}
            \pgfmathsetmacro{\yy}{\dS - (\row+1)*\dC}
            \fill[blockred] (\xx+0.01,\yy+0.01) rectangle (\xx+\dC-0.01,\yy+\dC-0.01);
          \fi\fi
        }
      }
      \node[anchor=north] at (\dS/2,-0.25) {$D_x$};
      \node[anchor=north, font=\footnotesize]
        at (\dS/2,-0.85) {banded, $N_x \times N_x$};
    \end{scope}

    \node[font=\bfseries, anchor=north]
      at (7.0,-1.6) {one batched line kernel: apply $D_x$ to all columns};
  \end{tikzpicture}}
  \caption{The multi-RHS reshape.  A three-dimensional $x$-sweep is
  exposed as one operation on a matrix of right-hand sides once the
  field is reshaped from $(N_x, N_y, N_z)$ to $(N_x, N_yN_z)$.  Each
  column of $U$ is one $x$-line.  Dense or element-local factors can be
  applied by GEMM; truly banded finite-difference factors should use a
  banded or stencil multi-RHS kernel rather than being densified.  The
  same reshape works for $y$- and $z$-sweeps after a transpose that
  puts the active direction into the leading axis.}
  \label{fig:blas3}
\end{figure}
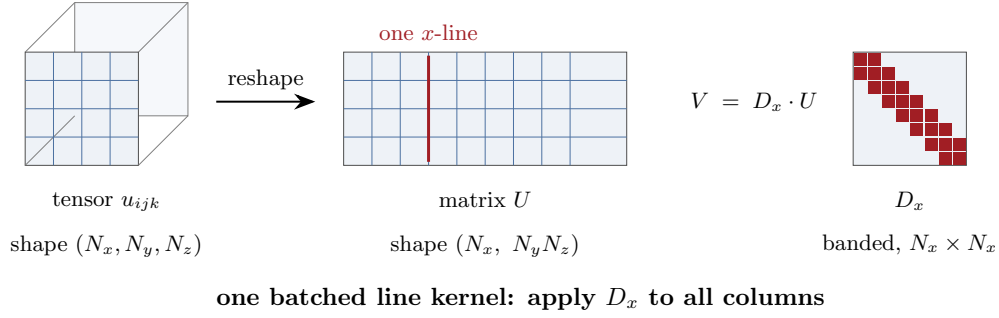

\Cref{fig:blas3} shows the production trick: reshape the tensor
$u \in \R^{N_x \times N_y \times N_z}$ into a matrix
$U \in \R^{N_x \times (N_y N_z)}$, whose columns are the $x$-lines.
The line apply is then the multi-right-hand-side operation
\begin{equation}
\label{eq:blas3}
  V \;=\; \Dx\, U.
\end{equation}
Equation~\eqref{eq:blas3} should be read with one important
implementation caveat.  If $\Dx$ is dense, spectral, or element-local
and small, this is exactly the shape of a BLAS-3 GEMM call.  If $\Dx$
is a narrow finite-difference stencil, densifying it would change the
work from $\bigO(w_x\Nall)$ to $\bigO(N_x\Nall)$ and destroy the
point of the method.  In that case, the correct production kernel is a
banded multi-RHS apply or a hand-written stencil kernel over the
columns of $U$.  Compact and Galerkin line solves use the solve analog:
factor the one-dimensional banded matrix once, then apply the factors
to many right-hand sides at once.  The same reshape works for the
$y$- and $z$-sweeps after a transpose or stride adjustment that places
the active direction first; see \cref{sec:memory} for the data-layout
side of the argument.

\paragraph{Why this matters.}
The reshape is the single most important reason production Cartesian
codes look fast.  It exposes reuse across many lines, lets dense or
element-local pieces fall into optimized GEMM, and lets banded pieces
use batched line kernels instead of scalar loops.  The
linear-algebraic content is unchanged; only the loop order and memory
layout are different.  Once the right kernel is chosen for the
bandwidth at hand, this reshape is essentially free.

\subsection{Trick 2: sum factorization (the spectral-element secret)}
\label{sec:sumfac}

A related, more dramatic speedup governs tensor-product Galerkin and
spectral-element methods at high polynomial order $p$.  The
right-hand-side quadrature sum on one element of $(p+1)^d$ points
in $d$ dimensions is
\begin{equation}
  F_{abc} \;=\; \sum_{q_1, q_2, q_3} w_{q_1 q_2 q_3}\,
  \phi_a(\xi_{q_1})\, \phi_b(\eta_{q_2})\, \phi_c(\zeta_{q_3})\,
  f(\xi_{q_1}, \eta_{q_2}, \zeta_{q_3}).
\end{equation}
If you read this literally, every one of the $(p+1)^3$ outputs sums
over $(p+1)^3$ inputs.  With $n=p+1$ points per direction, that is
$\bigO(n^{2d}) = \bigO(n^{6})$ operations per element in 3D.  For
$p=8$ ($n=9$), the naive count is about $5.3\times10^5$ multiply-adds
per element; for $p=16$ ($n=17$), it is about $2.4\times10^7$.  This
is the scaling that kept high-order methods out of production
engineering for two decades.

The cure lives in the same separation of variables that ran the rest
of the paper.  Do not sum all three indices at once; sum one at a
time.  Contract $q_3$ with $\phi_c$ to make an intermediate of shape
$(p{+}1)^2 \times (p{+}1)$; then contract $q_2$ with $\phi_b$; then
contract $q_1$ with $\phi_a$:
\begin{equation}
\label{eq:sumfac}
  F_{abc}
  \;=\;
  \sum_{q_1} \phi_a(\xi_{q_1})
  \underbrace{
  \sum_{q_2} \phi_b(\eta_{q_2})
  \underbrace{
    \sum_{q_3} \phi_c(\zeta_{q_3})\,
    w_{q_1 q_2 q_3}\, f_{q_1 q_2 q_3}
  }_{\text{cost } \bigO(n^{4})}
  }_{\text{cost } \bigO(n^{4})}
  \qquad\text{(and one more outer sum)}.
\end{equation}
Each contraction costs $(p{+}1)^3 \times (p{+}1) = \bigO(n^{4})$.
Three contractions, three directions.  The total drops from
$\bigO(n^{2d})$ to $\bigO(d\, n^{d+1})$, which is $\bigO(n^{4})$ in
3D.  Including the three contractions, the leading-count speedup is
about $n^2/3$: roughly $27\times$ at $p=8$ and $96\times$ at $p=16$
(\cref{fig:sumfac}).  This is why spectral-element and high-order DG
codes can run at orders that would otherwise be unthinkable~\citep{deville2002high}.
In Kronecker notation,
\eqref{eq:sumfac} is the contraction
$F = (\Phi_\xi \kron \Phi_\eta \kron \Phi_\zeta)^T (W\odot f)$,
evaluated as three sequential 1D matrix multiplies rather than as
one Kronecker-product matvec.  It is the same factorization that ran
the derivatives, the mass inversion, and the ADI substeps, applied
now to the quadrature loop itself.

\begin{figure}[t]
  \centering
  \begin{tikzpicture}[
    every node/.style={font=\small},
    box/.style={draw=black!70, fill=#1, minimum width=1.55cm, minimum height=0.9cm, inner sep=0pt, align=center},
    arr/.style={-{Stealth[length=3mm]}, line width=0.9pt, black!75},
  ]
    \definecolor{stage1}{HTML}{dbe7f5}
    \definecolor{stage2}{HTML}{e5d7eb}
    \definecolor{stage3}{HTML}{f8dfd9}
    \definecolor{stage4}{HTML}{d8efd8}
    \definecolor{bord}{HTML}{333333}

    \node[box=stage1] (S1) at (0,0) {$f_{q_1 q_2 q_3}$};
    \node[anchor=north, font=\footnotesize] at ([yshift=-0.05cm]S1.south) {all three $q$ indices};

    \node[box=stage2] (S2) at (3.7,0) {$G_{q_1 q_2 c}$};
    \node[anchor=north, font=\footnotesize] at ([yshift=-0.05cm]S2.south) {$q_3 \to c$};

    \node[box=stage3] (S3) at (7.4,0) {$H_{q_1 b c}$};
    \node[anchor=north, font=\footnotesize] at ([yshift=-0.05cm]S3.south) {$q_2 \to b$};

    \node[box=stage4] (S4) at (11.1,0) {$F_{abc}$};
    \node[anchor=north, font=\footnotesize] at ([yshift=-0.05cm]S4.south) {$q_1 \to a$};

    \draw[arr] (S1) -- node[above, font=\footnotesize]{\small$\sum_{q_3} \phi_c(\zeta_{q_3})\ldots$}
                         node[below, font=\footnotesize]{$\bigO(n^{4})$} (S2);
    \draw[arr] (S2) -- node[above, font=\footnotesize]{$\sum_{q_2} \phi_b(\eta_{q_2})\ldots$}
                         node[below, font=\footnotesize]{$\bigO(n^{4})$} (S3);
    \draw[arr] (S3) -- node[above, font=\footnotesize]{$\sum_{q_1} \phi_a(\xi_{q_1})\ldots$}
                         node[below, font=\footnotesize]{$\bigO(n^{4})$} (S4);

    \node[anchor=north, font=\small] at (5.55,-1.4)
      {total $\bigO(d\,n^{d+1}) = \bigO(n^{4})$ in 3D, versus $\bigO(n^{2d}) = \bigO(n^{6})$ done all at once};
  \end{tikzpicture}
  \caption{Sum factorization.  The $\bigO(n^{2d})$ nested quadrature
  sum of a spectral-element right-hand side, with $n=p+1$ points per
  direction, becomes three sequential contractions, each
  $\bigO(n^{d+1})$ in 3D.  For $p=16$ ($n=17$), the leading count is
  about $24$ million operations naively versus about $250$ thousand
  after three contractions.  The same dimension-by-dimension
  contraction that runs every sweep in this paper runs the quadrature,
  too.}
  \label{fig:sumfac}
\end{figure}

\subsection{Trick 3: pencil decomposition for parallel execution}
\label{sec:pencil}

\begin{figure}[t]
  \centering
  \resizebox{\linewidth}{!}{%
  \begin{tikzpicture}[
    x={(0.72cm,0cm)}, y={(0cm,0.72cm)}, z={(-0.48cm,-0.36cm)},
    every node/.style={font=\small},
  ]
    \pgfmathsetmacro{\Nx}{6}
    \pgfmathsetmacro{\Ny}{4}
    \pgfmathsetmacro{\Nz}{4}
    \definecolor{rankA}{HTML}{9ec5e8}
    \definecolor{rankB}{HTML}{ebb5b6}
    \definecolor{rankC}{HTML}{b9dbae}
    \definecolor{rankD}{HTML}{e6cf99}
    \definecolor{boxedge}{HTML}{555555}

    \newcommand{\wireframe}{%
      \draw[boxedge, line width=0.5pt] (0,0,0) -- (\Nx,0,0) -- (\Nx,\Ny,0) -- (0,\Ny,0) -- cycle;
      \draw[boxedge, line width=0.5pt] (0,0,\Nz) -- (\Nx,0,\Nz) -- (\Nx,\Ny,\Nz) -- (0,\Ny,\Nz) -- cycle;
      \draw[boxedge, line width=0.5pt] (0,0,0) -- (0,0,\Nz);
      \draw[boxedge, line width=0.5pt] (\Nx,0,0) -- (\Nx,0,\Nz);
      \draw[boxedge, line width=0.5pt] (\Nx,\Ny,0) -- (\Nx,\Ny,\Nz);
      \draw[boxedge, line width=0.5pt] (0,\Ny,0) -- (0,\Ny,\Nz);
    }

    \begin{scope}[xshift=0cm]
      \wireframe
      \foreach \jstart/\kstart/\col in {0/0/rankA, 2/0/rankB, 0/2/rankC, 2/2/rankD} {
        \fill[\col, opacity=0.85] (\Nx,\jstart,\kstart) -- (\Nx,\jstart+2,\kstart) -- (\Nx,\jstart+2,\kstart+2) -- (\Nx,\jstart,\kstart+2) -- cycle;
      }
      \foreach \jstart/\kstart/\col in {0/0/rankA, 2/0/rankB, 0/2/rankC, 2/2/rankD} {
        \draw[boxedge, line width=0.7pt] (\Nx,\jstart,\kstart) -- (\Nx,\jstart+2,\kstart) -- (\Nx,\jstart+2,\kstart+2) -- (\Nx,\jstart,\kstart+2) -- cycle;
      }
      \node[anchor=north] at (\Nx/2,-0.4,\Nz) {rank owns full $x$-line};
    \end{scope}

    \begin{scope}[xshift=7.5cm]
      \wireframe
      \foreach \istart/\kstart/\col in {0/0/rankA, 3/0/rankB, 0/2/rankC, 3/2/rankD} {
        \fill[\col, opacity=0.85] (\istart,\Ny,\kstart) -- (\istart+3,\Ny,\kstart) -- (\istart+3,\Ny,\kstart+2) -- (\istart,\Ny,\kstart+2) -- cycle;
      }
      \foreach \istart/\kstart/\col in {0/0/rankA, 3/0/rankB, 0/2/rankC, 3/2/rankD} {
        \draw[boxedge, line width=0.7pt] (\istart,\Ny,\kstart) -- (\istart+3,\Ny,\kstart) -- (\istart+3,\Ny,\kstart+2) -- (\istart,\Ny,\kstart+2) -- cycle;
      }
      \node[anchor=north] at (\Nx/2,-0.4,\Nz) {rank owns full $y$-line};
    \end{scope}

    \begin{scope}[xshift=14.0cm]
      \wireframe
      \foreach \istart/\jstart/\col in {0/0/rankA, 3/0/rankB, 0/2/rankC, 3/2/rankD} {
        \fill[\col, opacity=0.85] (\istart,\jstart,0) -- (\istart+3,\jstart,0) -- (\istart+3,\jstart+2,0) -- (\istart,\jstart+2,0) -- cycle;
      }
      \foreach \istart/\jstart/\col in {0/0/rankA, 3/0/rankB, 0/2/rankC, 3/2/rankD} {
        \draw[boxedge, line width=0.7pt] (\istart,\jstart,0) -- (\istart+3,\jstart,0) -- (\istart+3,\jstart+2,0) -- (\istart,\jstart+2,0) -- cycle;
      }
      \node[anchor=north] at (\Nx/2,-0.4,\Nz) {rank owns full $z$-line};
    \end{scope}

    \begin{scope}[x={(1cm,0cm)}, y={(0cm,1cm)}]
      \node[anchor=south, font=\bfseries\large] at (1.2,3.82) {$x$-pencils};
      \node[anchor=south, font=\bfseries\large] at (8.7,3.82) {$y$-pencils};
      \node[anchor=south, font=\bfseries\large] at (15.2,3.82) {$z$-pencils};

      \draw[line width=1.25pt, -{Stealth[length=3.2mm, width=2.7mm]}, black!70]
        (4.60,3.42) -- (6.10,3.42);
      \node[anchor=south, font=\footnotesize, fill=white, inner sep=1.6pt]
        at (5.35,3.50) {all-to-all};

      \draw[line width=1.25pt, -{Stealth[length=3.2mm, width=2.7mm]}, black!70]
        (12.10,3.42) -- (13.60,3.42);
      \node[anchor=south, font=\footnotesize, fill=white, inner sep=1.6pt]
        at (12.85,3.50) {all-to-all};
    \end{scope}
  \end{tikzpicture}}
  \caption{Pencil decomposition.  A $P_\alpha \times P_\beta$
  process grid owns the two perpendicular directions; each rank
  holds a complete ``pencil'' of the third direction.  Tiles in
  distinct colors represent distinct MPI ranks.  Between sweeps, a
  collective all-to-all transpose rotates the pencil orientation so
  that the next direction about to be swept is the local one.  The
  1D kernel inside a rank never changes; only which direction counts
  as local.}
  \label{fig:pencil}
\end{figure}
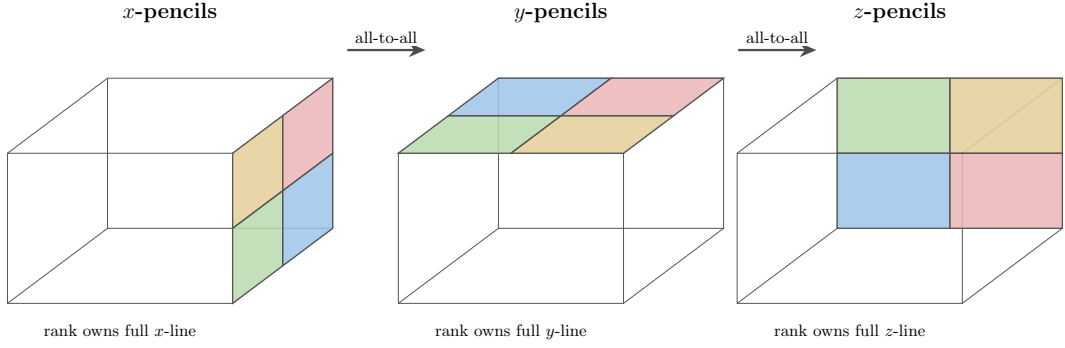

The serial story only takes us so far.  On a cluster of thousands
of MPI ranks, one direction of the grid will always be scattered
across ranks, and a sweep in that direction must move data.  The
cleanest response is the \emph{pencil decomposition}: lay the MPI
ranks on a two-dimensional process grid of size $P_\alpha \times
P_\beta$, so that the third direction lives entirely on each rank as
a long ``pencil'' of unknowns.  Before each sweep, a collective
transpose rotates the pencil orientation so that the direction about
to be swept is the local one (\cref{fig:pencil}).  A full
three-direction sweep therefore uses two all-to-all transposes; the
third orientation can often be reused for the next time step.

The FFT-based DNS solver
\texttt{CaNS}~\citep{costa2018fft} and the
\texttt{2DECOMP\&FFT} pencil-decomposition
library~\citep{li2010decomp,rolfo2023decomp}
both follow this pattern.  The pencil decomposition is the parallel
analog of the serial transpose that turns a $y$- or $z$-sweep into a
contiguous kernel (\cref{sec:memory}).  At both scales the idea is the
same: the tensor-product grid factors, and the algorithm is free to
factor with it.

\subsection{Three tricks, one idea}

The three tricks are one idea seen at three scales.  The multi-RHS
reshape exposes $N_y N_z$ line operations as one batched kernel because
the stride pattern of a sweep is separable in $(i,j,k)$.  The sum
factorization collapses $\bigO(n^{2d})$ integrals into
$\bigO(d\,n^{d+1})$ because the quadrature kernel is separable in
$(q_1,q_2,q_3)$.  The pencil decomposition collapses a distributed
sweep into a local one, with a transpose before and after, because the
data layout is separable across the MPI ranks.  Every speedup recovers
the same separable structure the PDE always had, at a scale the linear
algebra alone could not reach.

This is also the right design rule for software.  A Cartesian PDE code
should choose an active direction, expose all lines in that direction
as the columns of a matrix, apply the one-dimensional operator or solve
to all columns at once, and then rotate the layout for the next
direction.  The apparent bookkeeping around reshapes, transposes, and
MPI collectives is not auxiliary machinery; it is the mechanism that
lets the mathematical factorization survive contact with memory
hierarchy and distributed hardware.  Once this rule is followed, the
same source-level idea covers finite differences, compact schemes,
spline Galerkin methods, ADI steps, and FFT Poisson solvers.

\section{A numerical illustration: assembled versus matrix-free}
\label{sec:experiment}

The cost argument is easy to state and easy to check.  We solve the
Poisson problem $-\nabla^2 u = g$ on the unit cube with homogeneous
Dirichlet boundaries and the manufactured solution
$u = \sin(\pi x)\sin(\pi y)\sin(\pi z)$, discretized by the standard
seven-point stencil on an $N\times N\times N$ interior grid, and solve the
\emph{same} discrete system two ways.  The assembled route forms the
monolithic sparse Laplacian of size $N^3\times N^3$ as the Kronecker sum
\eqref{eq:3d-laplacian} and factors it with a sparse direct solver.  The
matrix-free route never builds that matrix: it solves by fast
diagonalization (\cref{sec:fast-diag})---a discrete sine transform along
each axis, a pointwise division by
$\lambda_i^x+\lambda_j^y+\lambda_k^z$, and an inverse transform.  Because
both routes solve the identical discrete system, any difference is in
storage and run time, not in the answer.

\begin{table}[t]
\centering
\footnotesize
\renewcommand{\arraystretch}{1.25}
\begin{tabular}{@{}rrrrrrr@{}}
\toprule
& & \multicolumn{2}{c}{assembled direct solve}
  & \multicolumn{2}{c}{matrix-free fast diag.} & \\
\cmidrule(lr){3-4}\cmidrule(lr){5-6}
$N$ & unknowns $N^3$ & factor (MB) & time (s)
    & storage (KB) & time (s) & $\|u_h-u\|_\infty$ \\
\midrule
$16$  & $4{,}096$       & $14.8$   & $0.08$  & $0.38$ & $0.0001$ & $2.8\times10^{-3}$ \\
$32$  & $32{,}768$      & $385$    & $3.45$  & $0.77$ & $0.0006$ & $7.5\times10^{-4}$ \\
$48$  & $110{,}592$     & $2{,}799$& $55.4$  & $1.15$ & $0.0022$ & $3.4\times10^{-4}$ \\
$64$  & $262{,}144$     & \multicolumn{2}{c}{infeasible} & $1.54$ & $0.0044$ & $2.0\times10^{-4}$ \\
$128$ & $2{,}097{,}152$ & \multicolumn{2}{c}{infeasible} & $3.07$ & $0.073$  & $4.9\times10^{-5}$ \\
$256$ & $16{,}777{,}216$& \multicolumn{2}{c}{infeasible} & $6.14$ & $1.98$   & $1.2\times10^{-5}$ \\
\bottomrule
\end{tabular}
\caption{The same 3D Poisson problem solved two ways.  ``Factor'' is the
fill-in of the sparse direct factorization of the assembled
$N^3\times N^3$ matrix; ``storage'' is the matrix-free operator, the three
sets of one-dimensional eigenvalues.  The error column is shared: both
routes return the identical discrete solution (they agree to about
$10^{-14}$), so $\|u_h-u\|_\infty$ is the same for each and falls at the
expected second-order rate.  The assembled factorization is run only while
it remains practical on a single workstation; past $48^3$ its fill-in
exhausts memory, while the matrix-free solver continues to
$256^3 \approx 1.7\times10^{7}$ unknowns.  Driver:
\texttt{scripts/poisson3d\_benchmark.py}.}
\label{tab:poisson3d}
\end{table}

\begin{figure}[t]
  \centering
  \includegraphics[width=\linewidth]{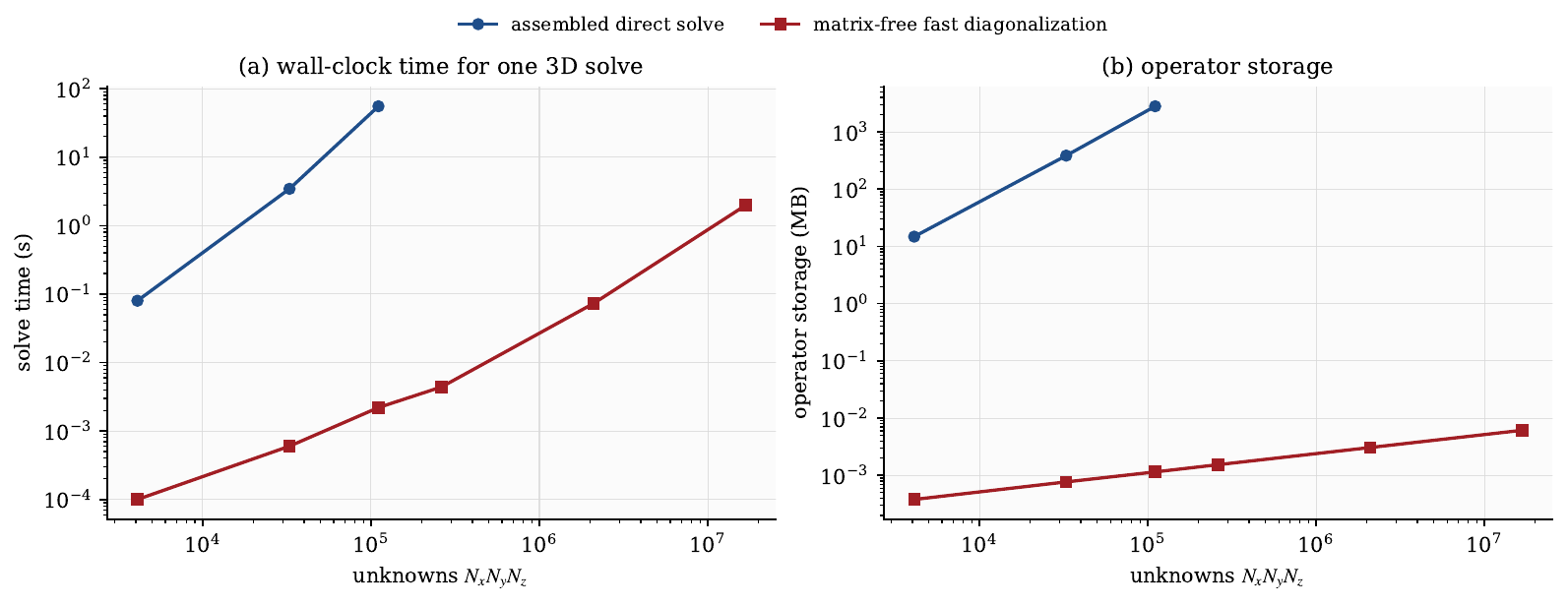}
  \caption{The same data as \cref{tab:poisson3d}, on log--log axes.  Never
  forming the three-dimensional matrix moves both axes by orders of
  magnitude: the matrix-free solver is faster (left) and stores a kilobyte
  of one-dimensional operators where the assembled factor stores gigabytes
  of fill-in (right).  The assembled curves stop where a direct
  factorization stops being practical; the matrix-free curves do not.}
  \label{fig:poisson3d}
\end{figure}

The numbers make three points.  First, the matrix-free route costs no
accuracy: the two solutions are the same discrete field to roughly
$10^{-14}$, and both converge at second order.  Second, the assembled
route pays for the three-dimensional matrix through fill-in.  At
$48^3 \approx 1.1\times10^{5}$ unknowns its sparse factor already needs
about $2.8$~GB, whereas the matrix-free operator is the set of
one-dimensional eigenvalues, about one kilobyte---a storage ratio above
two million.  Third, that fill-in is also a time wall: the direct
factorization takes nearly a minute at $48^3$ and exhausts a workstation
beyond it, while the matrix-free solver clears $256^3$, about seventeen
million unknowns, in roughly two seconds.  The $\bigO(\Nall)$ storage and
$\bigO(\Nall\log N_{\max})$ work claimed throughout this paper are not
asymptotic promises; they are visible at ordinary grid sizes.

\paragraph{Code availability.}
The scripts that regenerate every figure and \cref{tab:poisson3d} in this
paper, together with the \LaTeX{} source, are available at
\url{https://github.com/bay-yearick-lab/no-3d-matrices}.

\section{Summary and implementation recipe}
\label{sec:summary}

The practical conclusion of the paper fits on one page.  Let
$\Nall = N_x N_y N_z$ denote the total number of unknowns.  The
central result is this:
\emph{every local tensor-product derivative sweep and every
factorized tensor-product line solve costs $\bigO(\Nall)$ and stores
only the one-dimensional banded operators.}

\begin{enumerate}
  \item \textbf{Build 1D operators.}
    Construct $\Dx$, $\Dxx$ (and $\Ax$, $\Rx$ for compact schemes)
    in each direction.  For Galerkin or B-spline methods, build the 1D
    mass matrices $M_x$, $M_y$, $M_z$ and stiffness matrices $K_x$,
    $K_y$, $K_z$.  For spline collocation, build the banded nodal maps
    $B_x$, $B_x'$, $B_x''$ rather than relying on a formed product such
    as $B_x'' B_x^{-1}$.  These are all small banded matrices of sizes
    $N_x$, $N_y$, $N_z$.

  \item \textbf{Evaluate derivatives by sweeping lines.}
    To compute $\partial u/\partial x$: make the $x$-direction
    contiguous and apply $\Dx$ to every $x$-line.  The $y$- and
    $z$-directions follow by relabeling.  For compact schemes each
    line operation is a banded matvec plus a banded solve.  For
    Galerkin and spline discretizations, use the factorized banded
    maps along each line.  The full sweep costs $\bigO(\Nall)$.  For
    production performance, issue each sweep as one batched line kernel
    on the reshape $(N_\xi) \times (\Nall / N_\xi)$; use dense
    BLAS only when the one-dimensional factor is dense or element-local
    (\cref{sec:blas3}).

  \item \textbf{Implicit time stepping with ADI.}
    Split the implicit operator direction by direction.  Each substep
    is a batch of 1D banded solves, making the implicit method
    asymptotically the same cost as explicit, namely $\bigO(\Nall)$
    per time step.  The line factors are $I - \tfrac{\nu\dt}{2}\Dxx$
    for finite differences, $\Axx - \tfrac{\nu\dt}{2}\Rxx$ for
    compact schemes, and $M_x + \tfrac{\nu\dt}{2}K_x$ for Galerkin
    and spline methods.  Every one of them is banded.

  \item \textbf{Direct Poisson and Helmholtz solves.}
    Eigendecompose the 1D operators once at setup, transform
    direction by direction, solve pointwise, transform back.  When
    the 1D eigenvectors are trigonometric (uniform grid with
    Dirichlet, Neumann, or periodic boundaries), the transforms are
    FFTs and the total cost is $\bigO(\Nall \log N_{\max})$.  In
    general the transforms are small dense products per line.

  \item \textbf{Scale up.}
    For large problems, parallelize with a pencil decomposition
    (\cref{sec:pencil}) so every direction can be made contiguous in
    turn by a collective transpose.  For high-order Galerkin or
    spectral-element methods, evaluate every quadrature sum by sum
    factorization (\cref{sec:sumfac}) to replace the
    $\bigO(n^{2d})$ trap with the $\bigO(d\,n^{d+1})$ reality, where
    $n=p+1$ is the number of points per direction.
\end{enumerate}

For fixed stencil width or fixed polynomial degree, local sweeps and
factorized line solves cost $\bigO(\Nall)$ arithmetic and require only
$\bigO(N_x+N_y+N_z)$ operator storage, in addition to the
$\bigO(\Nall)$ storage for the fields themselves.  Direct separable
Poisson and Helmholtz solvers add transform costs, typically
$\bigO(\Nall\log N_{\max})$ with FFT-type bases.  The line operations
fan out across cores without coupling inside a sweep, fit naturally
into cache and accelerator memory hierarchies, and are exactly the
shape targeted by batched banded solvers and dense BLAS kernels when
the one-dimensional factors are dense or element-local.  Classical
textbooks arrive at related algorithms by different
roads~\citep{hirsch2007numerical,leveque2007finite,iserles2009first};
the route taken here is the Kronecker product, which trades a page of
index gymnastics for a single line of algebra and gives a sharp
baseline against which matrix-free, low-rank, and learned accelerators
can be judged.  For learned accelerators the test is twofold---cost
against this $\bigO(\Nall)$ baseline, and accuracy once the inputs drift
away from the training regime, where the distinction between fitting and
genuine generalization becomes decisive~\citep{bay2024generalization}.

\paragraph{A closing word.}
A separable three-dimensional Cartesian problem is really a
one-dimensional problem, dressed up for a Cartesian grid.  The
operator respects that product structure; the data layout, the line
kernel, and the parallel decomposition do too; and none of them ever
needs to see the monolithic three-dimensional matrix whose existence
caused the alarm in the first place.  When the Kronecker product,
the word ``line,'' and the phrase ``banded solve'' are held in the
same sentence, the picture is complete.  Every discretization family
this paper surveyed obeys it.  Long an unwritten habit of specialist
codes, it is elementary once the Kronecker product is taken as the
organizing principle.

\bibliographystyle{paper}
\bibliography{references}

\begin{thebibliography}{30}
\providecommand{\natexlab}[1]{#1}
\providecommand{\url}[1]{\texttt{#1}}
\expandafter\ifx\csname urlstyle\endcsname\relax
  \providecommand{\doi}[1]{doi: #1}\else
  \providecommand{\doi}{doi: \begingroup \urlstyle{rm}\Url}\fi

\bibitem[Li \& Laizet(2010)Li and Laizet]{li2010decomp}
Ning Li and Sylvain Laizet.
\newblock {2DECOMP\&FFT}: A highly scalable 2d decomposition library and {FFT}
  interface.
\newblock In \emph{Cray User Group 2010 Conference}, Edinburgh, United Kingdom,
  2010.

\bibitem[de~Boor(1978)]{deboor1978practical}
Carl de~Boor.
\newblock \emph{A Practical Guide to Splines}.
\newblock Springer, 1978.

\bibitem[Kim \& Moin(1985)Kim and Moin]{kim1985application}
John Kim and Parviz Moin.
\newblock Application of a fractional-step method to incompressible
  {N}avier--{S}tokes equations.
\newblock \emph{Journal of Computational Physics}, 59\penalty0 (2):\penalty0
  308--323, 1985.

\bibitem[Bay(2019)]{bay2019dissertation}
Yong~Yi Bay.
\newblock \emph{An Energy-Conservative Cut-Cell Method and Advanced
  {B}-Spline-Based Filtering Method for Flow Simulation}.
\newblock {Ph.D.} dissertation, University of Illinois at Urbana-Champaign,
  2019.
\newblock \url{https://hdl.handle.net/2142/106215}.

\bibitem[Bay \& Yearick(2026)Bay and Yearick]{bay2026spline}
Yong~Yi Bay and Kathleen~A. Yearick.
\newblock Solve for the hyperparameter, skip the search: {K}olmogorov-optimal
  scaling laws for spline regression.
\newblock \emph{arXiv preprint arXiv:2606.23575}, 2026.

\bibitem[Bay \& Yearick(2024)Bay and Yearick]{bay2024generalization}
Yong~Yi Bay and Kathleen~A. Yearick.
\newblock Machine learning vs deep learning: the generalization problem.
\newblock \emph{arXiv preprint arXiv:2403.01621}, 2024.

\bibitem[Van~Loan(2000)]{loan2000ubiquitous}
Charles~F Van~Loan.
\newblock The ubiquitous {K}ronecker product.
\newblock \emph{Journal of Computational and Applied Mathematics}, 123\penalty0
  (1--2):\penalty0 85--100, 2000.

\bibitem[Lele(1992)]{lele1992compact}
Sanjiva~K Lele.
\newblock Compact finite difference schemes with spectral-like resolution.
\newblock \emph{Journal of Computational Physics}, 103\penalty0 (1):\penalty0
  16--42, 1992.

\bibitem[Beam \& Warming(1976)Beam and Warming]{beam1976implicit}
Richard~M Beam and Robert~F Warming.
\newblock An implicit finite-difference algorithm for hyperbolic systems in
  conservation-law form.
\newblock \emph{Journal of Computational Physics}, 22\penalty0 (1):\penalty0
  87--110, 1976.

\bibitem[Golub \& Van~Loan(2013)Golub and Van~Loan]{golub2013matrix}
Gene~H Golub and Charles~F Van~Loan.
\newblock \emph{Matrix Computations}.
\newblock Johns Hopkins University Press, 4 edition, 2013.

\bibitem[Costa(2018)]{costa2018fft}
Pedro Costa.
\newblock A {FFT}-based finite-difference solver for massively-parallel direct
  numerical simulations of turbulent flows.
\newblock \emph{Computers \& Mathematics with Applications}, 76\penalty0
  (8):\penalty0 1853--1862, 2018.

\bibitem[Frigo \& Johnson(2005)Frigo and Johnson]{frigo2005design}
Matteo Frigo and Steven~G. Johnson.
\newblock The design and implementation of {FFTW3}.
\newblock \emph{Proceedings of the IEEE}, 93\penalty0 (2):\penalty0 216--231,
  2005.
\newblock \doi{10.1109/JPROC.2004.840301}.

\bibitem[Kolda \& Bader(2009)Kolda and Bader]{kolda2009tensor}
Tamara~G Kolda and Brett~W Bader.
\newblock Tensor decompositions and applications.
\newblock \emph{{SIAM} Review}, 51\penalty0 (3):\penalty0 455--500, 2009.

\bibitem[Lynch et~al.(1964)Lynch, Rice, and Thomas]{lynch1964direct}
Robert~E Lynch, John~R Rice, and David~H Thomas.
\newblock Direct solution of partial difference equations by tensor product
  methods.
\newblock \emph{Numerische Mathematik}, 6\penalty0 (1):\penalty0 185--199,
  1964.

\bibitem[Rolfo et~al.(2023)Rolfo, Flageul, Bartholomew, Spiga, and
  Laizet]{rolfo2023decomp}
Stefano Rolfo, C{\'e}dric Flageul, Paul Bartholomew, Filippo Spiga, and Sylvain
  Laizet.
\newblock The {2DECOMP\&FFT} library: an update with new {CPU/GPU}
  capabilities.
\newblock \emph{Journal of Open Source Software}, 8\penalty0 (91):\penalty0
  5813, 2023.
\newblock \doi{10.21105/joss.05813}.

\bibitem[Bewley et~al.(2001)Bewley, Moin, and Temam]{bewley2001dns}
Thomas~R Bewley, Parviz Moin, and Roger Temam.
\newblock {DNS}-based predictive control of turbulence: an optimal benchmark
  for feedback algorithms.
\newblock \emph{Journal of Fluid Mechanics}, 447:\penalty0 179--225, 2001.

\bibitem[Hughes et~al.(2005)Hughes, Cottrell, and
  Bazilevs]{hughes2005isogeometric}
Thomas~JR Hughes, John~A Cottrell, and Yuri Bazilevs.
\newblock Isogeometric analysis: {CAD}, finite elements, {NURBS}, exact
  geometry and mesh refinement.
\newblock \emph{Computer Methods in Applied Mechanics and Engineering},
  194\penalty0 (39--41):\penalty0 4135--4195, 2005.

\bibitem[Hirsch(2007)]{hirsch2007numerical}
Charles Hirsch.
\newblock \emph{Numerical Computation of Internal and External Flows}.
\newblock Butterworth-Heinemann, 2 edition, 2007.

\bibitem[Strang(2007)]{strang2007computational}
Gilbert Strang.
\newblock \emph{Computational Science and Engineering}.
\newblock Wellesley-Cambridge Press, 2007.

\bibitem[Thomas(1949)]{thomas1949elliptic}
Llewellyn~H Thomas.
\newblock Elliptic problems in linear difference equations over a network.
\newblock \emph{Watson Scientific Computing Laboratory Report}, 1949.

\bibitem[Douglas(1962)]{douglas1962alternating}
Jim Douglas.
\newblock Alternating direction methods for three space variables.
\newblock \emph{Numerische Mathematik}, 4\penalty0 (1):\penalty0 41--63, 1962.

\bibitem[Douglas \& Gunn(1964)Douglas and Gunn]{douglas1964numerical}
Jim Douglas, Jr and James~E Gunn.
\newblock A general formulation of alternating direction methods.
\newblock \emph{Numerische Mathematik}, 6\penalty0 (1):\penalty0 428--453,
  1964.

\bibitem[Iserles(2009)]{iserles2009first}
Arieh Iserles.
\newblock \emph{A First Course in the Numerical Analysis of Differential
  Equations}.
\newblock Cambridge University Press, 2 edition, 2009.

\bibitem[Deville et~al.(2002)Deville, Fischer, and Mund]{deville2002high}
Michel~O Deville, Paul~F Fischer, and Ernest~H Mund.
\newblock \emph{High-Order Methods for Incompressible Fluid Flow}.
\newblock Cambridge University Press, 2002.

\bibitem[Johnson(2005)]{johnson2005higher}
Richard~W Johnson.
\newblock Higher order {B}-spline collocation at the {G}reville abscissae.
\newblock \emph{Applied Numerical Mathematics}, 52\penalty0 (1):\penalty0
  63--75, 2005.

\bibitem[LeVeque(2007)]{leveque2007finite}
Randall~J LeVeque.
\newblock \emph{Finite Difference Methods for Ordinary and Partial Differential
  Equations}.
\newblock SIAM, 2007.

\bibitem[Cottrell et~al.(2009)Cottrell, Hughes, and
  Bazilevs]{cottrell2009isogeometric}
J~Austin Cottrell, Thomas~JR Hughes, and Yuri Bazilevs.
\newblock \emph{Isogeometric Analysis: Toward Integration of {CAD} and {FEA}}.
\newblock Wiley, 2009.

\bibitem[Fornberg(1988)]{fornberg1988generation}
Bengt Fornberg.
\newblock Generation of finite difference formulas on arbitrarily spaced grids.
\newblock \emph{Mathematics of Computation}, 51\penalty0 (184):\penalty0
  699--706, 1988.

\bibitem[Peaceman \& Rachford(1955)Peaceman and
  Rachford]{peaceman1955numerical}
Donald~W Peaceman and Henry~H Rachford, Jr.
\newblock The numerical solution of parabolic and elliptic differential
  equations.
\newblock \emph{Journal of the Society for Industrial and Applied Mathematics},
  3\penalty0 (1):\penalty0 28--41, 1955.

\bibitem[Haidvogel \& Zang(1979)Haidvogel and Zang]{haidvogel1979accurate}
Dale~B Haidvogel and Thomas Zang.
\newblock The accurate solution of {P}oisson's equation by expansion in
  {C}hebyshev polynomials.
\newblock \emph{Journal of Computational Physics}, 30\penalty0 (2):\penalty0
  167--180, 1979.

\end{thebibliography}

\clearpage
\begin{appendices}

\section{Algebraic details for the 1D reductions}
\label{app:linalg}

The main text uses the Kronecker product as a language for loops.
This appendix gives the short, self-contained derivations that
translate that language back into arrays, line solves, and direct
diagonalization formulas.  The whole edifice rests on two identities:
the mixed-product rule of \cref{app:mixed-product}, and the
\emph{vec identity} of \cref{app:vec-identity}.  Every subsequent
derivation in this appendix is one or two applications of one or both.

\subsection{Mixed-product rule and Kronecker inverse}
\label{app:mixed-product}
\label{app:kron-inverse}

\paragraph{Claim.}
For matrices $A \in \R^{m \times n}$, $C \in \R^{n \times \ell}$,
$B \in \R^{p \times q}$, $D \in \R^{q \times r}$,
\begin{equation}
\label{eq:mixed-product-app}
  (A \kron B)(C \kron D) \;=\; (AC) \kron (BD).
\end{equation}

\paragraph{Proof.}
By block matrix multiplication, the $(i,k)$ block of a product of two
block matrices is $\sum_{j} (\text{block}_{ij}^{(1)})\,(\text{block}_{jk}^{(2)})$.
By definition of the Kronecker product, the $(i,j)$ block of
$A \kron B$ is the scalar-times-matrix $a_{ij}\,B$, and similarly the
$(j,k)$ block of $C \kron D$ is $c_{jk}\,D$.  The $(i,k)$ block of the
left-hand side of \eqref{eq:mixed-product-app} is therefore
\begin{equation*}
  \sum_{j=1}^{n} (a_{ij}\,B)(c_{jk}\,D)
  \;=\;
  \sum_{j=1}^{n} a_{ij}\, c_{jk}\, B\, D
  \;=\;
  \biggl(\sum_{j=1}^{n} a_{ij}\,c_{jk}\biggr) B\, D
  \;=\;
  (AC)_{ik}\,(BD),
\end{equation*}
where the second equality uses the fact that $a_{ij}$ and $c_{jk}$
are scalars and slide through $B$ and $D$ freely.  The right-hand
side is the $(i,k)$ block of $(AC) \kron (BD)$, which establishes
\eqref{eq:mixed-product-app}.  $\square$

\paragraph{Corollary (inverse rule).}
If $A$ and $B$ are square and invertible, take $C = A^{-1}$ and
$D = B^{-1}$ in \eqref{eq:mixed-product-app}:
\begin{equation*}
  (A \kron B)(A^{-1} \kron B^{-1})
  \;=\; (A A^{-1}) \kron (B B^{-1})
  \;=\; I_m \kron I_p
  \;=\; I_{mp},
\end{equation*}
and the product taken in the other order yields $I_{mp}$ as well, so
the inverse is two-sided.  Hence
\begin{equation}
\label{eq:kron-inv-app}
  (A \kron B)^{-1} \;=\; A^{-1} \kron B^{-1}.
\end{equation}
The triple-product extension follows by associativity:
$(A \kron B \kron C)^{-1} = A^{-1} \kron B^{-1} \kron C^{-1}$
whenever each factor is square and invertible.

\subsection{The vec identity (master lemma)}
\label{app:vec-identity}

Throughout the paper, $\vect$ uses the column-major convention.  For
$U \in \R^{m \times n}$ with columns $U_{:,1}, U_{:,2}, \ldots, U_{:,n}$,
\begin{equation}
\label{eq:vec-def}
  \vect(U)
  \;=\;
  \begin{bmatrix}
    U_{:,1} \\[2pt]
    U_{:,2} \\[2pt]
    \vdots \\[2pt]
    U_{:,n}
  \end{bmatrix}
  \;\in\; \R^{mn}.
\end{equation}
The columns of $U$ stack vertically into one tall column vector.
Reading $\vect(U)$ from top to bottom walks the row index ($x$)
fastest and the column index ($y$) slowest, which is the ordering
used throughout the main text.

\paragraph{Claim (vec identity).}
For $A \in \R^{m \times m}$, $B \in \R^{n \times n}$, and
$X \in \R^{m \times n}$,
\begin{equation}
\label{eq:vec-identity}
  \boxed{\,(B^T \kron A)\,\vect(X) \;=\; \vect(A X B).\,}
\end{equation}

\paragraph{Proof.}
Write $Y = A X B$.  Using $B_{:,j} = \sum_k B_{kj}\,e_k$ where $e_k$
is the $k$th standard basis vector and $X e_k = X_{:,k}$, the $j$th
column of $Y$ is
\begin{equation*}
  Y_{:,j}
  \;=\; A\, X\, B_{:,j}
  \;=\; A\,\sum_{k=1}^{n} B_{kj}\, X_{:,k}
  \;=\; \sum_{k=1}^{n} B_{kj}\,(A\, X_{:,k}).
\end{equation*}
Stacking the columns of $Y$ in column-major order gives
\begin{equation*}
  \vect(Y)
  \;=\;
  \begin{bmatrix}
    \sum_{k} B_{k 1}\, A\, X_{:,k} \\[2pt]
    \sum_{k} B_{k 2}\, A\, X_{:,k} \\[2pt]
    \vdots \\[2pt]
    \sum_{k} B_{k n}\, A\, X_{:,k}
  \end{bmatrix}.
\end{equation*}
This vector is exactly $(B^T \kron A)\,\vect(X)$ read by blocks: the
$j$th block is
$\sum_k (B^T)_{j k}\, A\, X_{:,k} = \sum_k B_{k j}\, A\, X_{:,k}$,
as required.  $\square$

Two specializations of \eqref{eq:vec-identity} carry most of the
computational weight in the main text.  Setting $B = I_n$ kills the
right multiplication and gives
\begin{equation}
\label{eq:vec-IA}
  (I_n \kron A)\,\vect(X) \;=\; \vect(A X).
\end{equation}
Setting $A = I_m$ instead, and writing $\tilde B = B^T$, gives
\begin{equation}
\label{eq:vec-BI}
  (\tilde B \kron I_m)\,\vect(X) \;=\; \vect(X\, \tilde B^T).
\end{equation}
In words: \eqref{eq:vec-IA} says a matrix in the rightmost Kronecker
slot acts down the columns of the array form of $\vect(X)$;
\eqref{eq:vec-BI} says a matrix in the next slot acts across the
rows.  These two statements are the entire algebraic content of the
2D line-sweep picture.

\subsection{The 3D directional actions}
\label{app:vec}

In three dimensions, store the field as $U(:,:,:) \in \R^{N_x \times
N_y \times N_z}$ with the $x$-index first.  Write $U^{(k)} = U(:,:,k)
\in \R^{N_x \times N_y}$ for the $k$th $x$-$y$ plane.  The flattened
vector $\bm{u} \in \R^{\Nall}$ is built by stacking the columns
within each plane, then stacking the planes:
\begin{equation}
\label{eq:3d-flat}
  \bm{u}
  \;=\;
  \begin{bmatrix}
    \vect(U^{(1)}) \\[2pt]
    \vect(U^{(2)}) \\[2pt]
    \vdots \\[2pt]
    \vect(U^{(N_z)})
  \end{bmatrix}
  \;\in\; \R^{\Nall}.
\end{equation}
With this convention, $x$ varies fastest, then $y$, then $z$.  The
three directional actions of the main text follow immediately from
the 2D specializations \eqref{eq:vec-IA}--\eqref{eq:vec-BI} applied
plane by plane.

\paragraph{$x$-direction.}
The outer factor $\Iz$ leaves the plane index $k$ untouched, so the
$k$th block of $(\Iz \kron \Iy \kron \Dx)\,\bm{u}$ is
$(\Iy \kron \Dx)\,\vect(U^{(k)})$.  By \eqref{eq:vec-IA} this equals
$\vect(\Dx\,U^{(k)})$, which gives
\begin{equation}
\label{eq:app-x-action}
  (\Iz \kron \Iy \kron \Dx)\,\bm{u}
  \;=\;
  \begin{bmatrix}
    \vect(\Dx\,U^{(1)}) \\[2pt]
    \vdots \\[2pt]
    \vect(\Dx\,U^{(N_z)})
  \end{bmatrix}.
\end{equation}
The interpretation is the line sweep: apply $\Dx$ to every column of
every plane.

\paragraph{$y$-direction.}
The $k$th block of $(\Iz \kron \Dy \kron \Ix)\,\bm{u}$ is
$(\Dy \kron \Ix)\,\vect(U^{(k)})$.  By \eqref{eq:vec-BI} with
$\tilde B = \Dy$ this equals $\vect(U^{(k)}\,\Dy^T)$, so
\begin{equation}
\label{eq:app-y-action}
  (\Iz \kron \Dy \kron \Ix)\,\bm{u}
  \;=\;
  \begin{bmatrix}
    \vect(U^{(1)}\,\Dy^T) \\[2pt]
    \vdots \\[2pt]
    \vect(U^{(N_z)}\,\Dy^T)
  \end{bmatrix}.
\end{equation}
Multiplying $U^{(k)}$ by $\Dy^T$ on the right is the same as applying
$\Dy$ along each row of $U^{(k)}$, that is, along each $y$-line at
fixed $z = k$.

\paragraph{$z$-direction.}
The $\Dz$ factor now sits on the plane index, so it produces a linear
combination of whole planes:
\begin{equation}
\label{eq:app-z-action}
  (\Dz \kron \Iy \kron \Ix)\,\bm{u}
  \;=\;
  \begin{bmatrix}
    \sum_{r=1}^{N_z} (\Dz)_{1 r}\,\vect(U^{(r)}) \\[2pt]
    \vdots \\[2pt]
    \sum_{r=1}^{N_z} (\Dz)_{N_z\, r}\,\vect(U^{(r)})
  \end{bmatrix}.
\end{equation}
Equivalently, treat the matrix
$M = [\,\vect(U^{(1)}) \mid \vect(U^{(2)}) \mid \cdots \mid \vect(U^{(N_z)})\,]
\in \R^{(N_x N_y) \times N_z}$,
whose $k$th column holds the flattened plane $U^{(k)}$.  Then
\eqref{eq:app-z-action} is the matrix product $M\,\Dz^T$ read column
by column, which is one banded matvec per $z$-line of length $N_z$.

These three formulas, \eqref{eq:app-x-action}--\eqref{eq:app-z-action},
are exactly the three sweeps of \cref{fig:hero-sweep} and the
\texttt{sweep} pseudocode of \cref{fig:pseudocode}, written without
indices.

\subsection{Kronecker sums and fast diagonalization}
\label{app:fast-diag-derivation}

The discrete Laplacian is the Kronecker sum
\begin{equation}
\label{eq:app-laplacian-kron-sum}
  L
  \;=\; \Iz \kron \Iy \kron \Dxx
  \;+\; \Iz \kron \Dyy \kron \Ix
  \;+\; \Dzz \kron \Iy \kron \Ix.
\end{equation}
Suppose each 1D second-derivative matrix is diagonalizable as
$\Dxx = S_x \Lambda_x S_x^{-1}$, $\Dyy = S_y \Lambda_y S_y^{-1}$, and
$\Dzz = S_z \Lambda_z S_z^{-1}$.  Each Kronecker term factors by the
mixed-product rule applied twice.  For the first term, write each
identity as $S_\xi S_\xi^{-1}$:
\begin{align}
  \Iz \kron \Iy \kron \Dxx
  &\;=\;
  (S_z S_z^{-1}) \kron (S_y S_y^{-1}) \kron (S_x \Lambda_x S_x^{-1})
  \nonumber \\
  &\;=\;
  (S_z \kron S_y \kron S_x)\,
  (\Iz \kron \Iy \kron \Lambda_x)\,
  (S_z^{-1} \kron S_y^{-1} \kron S_x^{-1}),
  \label{eq:fastdiag-step}
\end{align}
where the second line is one application of
\eqref{eq:mixed-product-app} extended to triple products by
associativity (the binary rule applies to any factorization
$XYZ = (XY)Z = X(YZ)$, so it iterates without change).  The other two
Kronecker terms in \eqref{eq:app-laplacian-kron-sum} factor identically
with $\Lambda_y$ and $\Lambda_z$ in the appropriate slot.  Crucially,
all three terms share the \emph{same} outer factors
$S_z \kron S_y \kron S_x$ and
$S_z^{-1} \kron S_y^{-1} \kron S_x^{-1}$, so the sum collects:
\begin{equation}
\label{eq:fastdiag-factored}
  L
  \;=\;
  (S_z \kron S_y \kron S_x)\,
  \underbrace{\bigl(\Iz \kron \Iy \kron \Lambda_x \;+\; \Iz \kron \Lambda_y \kron \Ix \;+\; \Lambda_z \kron \Iy \kron \Ix\bigr)}_{\text{diagonal}}\,
  (S_z^{-1} \kron S_y^{-1} \kron S_x^{-1}).
\end{equation}
The middle factor is diagonal because each summand is a Kronecker
product of diagonal matrices and is therefore diagonal, and the sum
of diagonal matrices is diagonal.  Its $(i,j,k)$ entry is
$\lambda_i^x + \lambda_j^y + \lambda_k^z$.

The Poisson solve $L\,\bm{u} = \bm{f}$ therefore reduces, in the
tensor-product eigenbasis
$\hat{\bm{u}} = (S_z^{-1} \kron S_y^{-1} \kron S_x^{-1})\,\bm{u}$,
to the scalar division
\begin{equation}
\label{eq:app-fastdiag-final}
  \hat u_{ijk}
  \;=\;
  \frac{\hat f_{ijk}}{\lambda_i^x + \lambda_j^y + \lambda_k^z},
\end{equation}
provided the denominator never vanishes, or else the corresponding
compatibility condition is imposed and the null mode is fixed by a
normalization.  A shifted operator $\alpha I + \beta L$ replaces the
denominator with
$\alpha + \beta(\lambda_i^x + \lambda_j^y + \lambda_k^z)$ and is
otherwise identical.

\subsection{Compact and spline formulas as line systems}
\label{app:factorized-lines}

The compact $x$-derivative system in 3D, written for the vector of
derivative values $\bm{v}$, is
\begin{equation*}
  (\Iz \kron \Iy \kron \Ax)\,\bm{v}
  \;=\;
  (\Iz \kron \Iy \kron \Rx)\,\bm{u}.
\end{equation*}
Read this plane by plane.  The outer $\Iz$ leaves the plane index
untouched; on each plane, the action of $\Iy \kron \Ax$ is
\eqref{eq:vec-IA} applied with $A = \Ax$ and $X = U^{(k)}$, and
similarly for $\Iy \kron \Rx$.  Thus for every plane index $k$,
\begin{equation}
\label{eq:app-compact-plane}
  \Ax\, V^{(k)} \;=\; \Rx\, U^{(k)}.
\end{equation}
Equation \eqref{eq:app-compact-plane} is itself a matrix equation;
reading it column by column (each column is one $x$-line at fixed
$y$ and $z$) gives the independent line systems
\begin{equation}
\label{eq:app-compact-line}
  \Ax\, v(:,j,k)
  \;=\;
  \Rx\, u(:,j,k),
  \qquad
  j = 1, \ldots, N_y,\;
  k = 1, \ldots, N_z.
\end{equation}
Each line is one banded matvec by $\Rx$ followed by one tridiagonal
solve with $\Ax$, executed by the Thomas algorithm of
\eqref{eq:thomas-fwd}--\eqref{eq:thomas-bwd}.  The $y$- and
$z$-direction compact derivatives factor identically, using
\eqref{eq:vec-BI} and the $z$-plane action of
\eqref{eq:app-z-action}, respectively.

The spline coefficient-recovery and weak-form derivative equations
have the same structure.  The 3D collocation system
$(\Iz \kron \Iy \kron B_x)\,\bm{\alpha} = \bm{u}$ becomes, by
\eqref{eq:vec-IA} on each plane,
\begin{equation}
  B_x\, A^{(k)} \;=\; U^{(k)},
  \qquad\text{or column-wise}\qquad
  B_x\, \alpha(:,j,k) \;=\; u(:,j,k),
\end{equation}
where $\alpha(:,j,k)$ is the column of coefficients along the
$x$-line at $(j,k)$.  Once recovered, the derivative samples on the
same line are
\begin{equation*}
  u_x(:,j,k) \;=\; B_x'\,\alpha(:,j,k).
\end{equation*}
The IGA coefficient-derivative system
\begin{equation*}
  (\Iz \kron \Iy \kron M_x)\,\bm{\beta}
  \;=\;
  (\Iz \kron \Iy \kron G_x)\,\bm{\alpha}
\end{equation*}
becomes, in the same way,
\begin{equation}
  M_x\, \beta(:,j,k)
  \;=\;
  G_x\, \alpha(:,j,k).
\end{equation}
In every case the multidimensional structure expands the number of
independent lines, not the size or shape of the kernel.  $\square$

\end{appendices}

\end{document}